\documentclass[11pt,twoside,epsf]{article}
\usepackage{amsmath,amssymb,amsthm}
\usepackage{bm}
\usepackage{amssymb}
\usepackage{color}
\usepackage{amsthm}
\usepackage{galois}
\usepackage{subfigure}
\usepackage{graphicx}

\newfont{\bb}{msbm10}

\newtheorem{remark}{Remark}[section]
\newtheorem{theorem}{Theorem}[section]

\makeatletter % '@' is now a normal "letter" for TeX

\@addtoreset{equation}{section}
\makeatother % '@' is restored as a "non-letter" character for TeX
\begin{document}
\cleardoublepage
\pagestyle{myheadings}

\bibliographystyle{plain}

\title{Explicit K-symplectic methods for nonseparable non-canonical Hamiltonian systems}
\author
{Beibei Zhu \\
{\it Department of Applied Mathematics, School of Mathematics and Physics}\\
{\it University of Science and Technology Beijing, Beijing 100083, China}\\
{\it Email: zhubeibei@lsec.cc.ac.cn}\\[2mm]
Lun Ji \\
{\it LSEC, ICMSEC, Academy of Mathematics and Systems Science,}\\
{\it Chinese Academy of Sciences, Beijing 100190, China}\\
{\it School of Mathematical Sciences, University of Chinese Academy of Sciences}\\
{\it  Beijing 100049, China}\\
{\it Email: ujeybn@lsec.cc.ac.cn}
 \\[2mm]
Aiqing Zhu \\
{\it LSEC, ICMSEC, Academy of Mathematics and Systems Science,}\\
{\it Chinese Academy of Sciences, Beijing 100190, China}\\
{\it School of Mathematical Sciences, University of Chinese Academy of Sciences}\\
{\it Beijing 100049, China}\\
{\it Email: zaq@lsec.cc.ac.cn}
\\[2mm]
Yifa Tang\\
{\it LSEC, ICMSEC, Academy of Mathematics and Systems Science}\\
{\it Chinese Academy of Sciences, Beijing 100190, China}\\
{\it School of Mathematical Sciences, University of Chinese Academy of Sciences}\\
{\it Beijing 100049, China} \\
{\it Email: tyf@lsec.cc.ac.cn}}

\date{}

\maketitle

\markboth{\small B.B. Zhu, L. Ji, A.Q. Zhu and Y.F. Tang}
{\small Poisson Integrators for Poisson systems}

\begin{abstract}
We propose efficient numerical methods for nonseparable non-canonical Hamiltonian systems which are explicit, K-symplectic in the extended phase space with long time energy conservation properties. They are based on extending the original phase space to several copies of the phase space and imposing a mechanical restraint on the copies of the phase space. Explicit K-symplectic methods are constructed for three non-canonical Hamiltonian systems. Numerical results show that they outperform the higher order Runge-Kutta methods in preserving the phase orbit and the energy of the system over long time.

\bigskip


\bigskip

\end{abstract}

\section{Introduction}

We propose efficient numerical methods for nonseparable non-canonical Hamiltonian systems which are explicit,
K-symplectic in the extended phase space with long time energy conservation properties. The idea is to extend the original phase space to several copies of the phase space and impose a mechanical restraint on the copies of the phase space. Explicit K-symplectic methods are constructed for three non-canonical Hamiltonian systems. Numerical results show that they outperform the higher order Runge-Kutta methods in preserving the phase orbit and the energy of the system over long time.

Many mechanical systems can be expressed as non-canonical Hamiltonian systems, such as the Lotka-Volterra model\cite{GNT}, the nonlinear Schr\"{o}dinger equation\cite{Brugnano,Tang1,Tang2}, the charged particle system\cite{He,Zhang2}, the guiding center system\cite{Littlejohn,Qin,Zhang,Zhu3}, the Maxwell-Vlasov equations\cite{Marsden,Morrison} and the ideal MHD formualtion\cite{Morrison1}. As the non-canonical Hamiltonian systems are generalizations of canonical Hamiltonian systems\cite{Arnold,FPS,SSC94}. They have non-canonical symplectic structure which are preserved by K-symplectic methods. K-symplectic methods exhibit advantageous energy preservation properties just like the symplectic methods for Hamiltonian systems\cite{Channell,Feng1984,FQ,Forest,GNT,Tang1}. Therefore, K-symplectic methods are preferred methods for long time simulations of non-canonical Hamiltonian systems. One approach to constructing K-symplectic methods for a non-canonical system is to transform the system to a canonical one by a coordinate transformtion\cite{Arnold,GNT} and use the symplectic method. The other approach is the generating function method\cite{FPS}, but it is also inevitable to seek for one coordinate transformation. To avoid the difficulty of finding the coordinate transformation, we prefer to use the splitting method for non-canonical Hamiltonian systems.

The splitting method is widely used for separable Hamiltonian systems\cite{Blanes,GNT,McL2} and non-canonical Hamiltonian systems\cite{He,Zhu,Zhu2} to construct numerical methods. Blanes and Moan used the splitting method to construct an efficient fourth order symplectic method for non-autonomous Hamiltonian systems\cite{Blanes}. He et al. used the splitting method to construct the explicit K-sympletic methods for the non-canonical charged particle system. Zhu et al. constructed the explicit K-symplectic methods for both the bright solitons motion and the dark solitons motion of the nonlinear Schr\"{o}dinger equation\cite{Zhu,Zhu2}. The explicit symplectic methods and K-symplectic methods have been constructed for separable systems while the work is much less for nonseparable systems. For several subclasses of nonseparable Hamiltonian systems, Pihajoki constructed the explicit symplectic methods based on the idea of extending the phase space by making two copies of it and using the splitting method\cite{Pihajoki}. To overcome the short-term simulation problem of Pihajoki's work, Tao imposed a mechanical restraint on the two copies of the phase space and successfully realized the long-term simulation\cite{Tao}. Our work is based on the pioneer work of Pihajoki and Tao, we construct the explicit K-symplectic methods for nonseparable non-canonical Hamiltonian systems by using the splitting method. As the non-canonical systems are more complicate than canonical systems, in some cases the two copies of the phase space are not enough to construct explicit K-symplectic methods, thus we need to make more copies.

Nonseparable systems are important gradients of non-canonical Hamiltonian systems. The guiding center system\cite{Littlejohn,Qin,Zhang}, the reduced magnetohydrodynamics equations\cite{Kaneko} and the Ablowitz-Ladik model of the nonlinear Schr\"{o}dinger equation under periodic condition\cite{GNT} are all nonseparable non-canonical Hamiltonian systems. Thus, it is critical to construct K-symplectic methods with long-term conservation property for nonseparable non-canonical systems. The idea is to extend the phase space and make several copies of the phase space and impose an artificial restraint for the several copies of the phase space. The constructed method is explicit and K-symplectic in the extended phase space. We extend the non-canonical systems in two ways according to their structures and analyze in which situation that the explicit K-symplectic methods can be constructed. To show their superiority in structure preservation, they are compared with higher order explicit Runge-Kutta methods\cite{Butcher1,Butcher2}.

The paper is organized as follows. Section 2 is focus on constructing K-symplectic methods for nonseparable systems with special structure and analyzing the situation in which the explicit K-symplectic methods can be constructed. In Section 3, the nonseparable systems are general systems and we make several copies of the phase space in order to construct explicit K-symplectic methods. In Section 4, numerical results in three nonseparable systems are provided. Finally, we give a short conclusion in Section 5.

\section{Explicit K-symplectic methods for nonseparable systems with special structure}
\label{Kspecial}

Consider the non-canonical Hamiltonian system
\begin{equation}\label{noncanHam}
\dot{Z}=K^{-1}(Z)\nabla H(Z),\quad z=(P,Q)^{\top}\in \mathbb{R}^{2n}
\end{equation}
where the skew-symmetric matrix satisfies
\begin{equation}\label{Kstructure}
K^{-1}(Z)=\bordermatrix{&\cr
 &O_n& -K_{12}(P,Q)\cr
 &K_{12}^{\top}(P,Q)&O_n\cr
 },
\end{equation}
we approximate its flow of an arbitrary nonseparable $H(P,Q)$. We aim at constructing higher order explicit K-symplectic methods by considering an augmented Hamiltonian\cite{Tao}
$$
\bar{H}(p,q,x,y):=H_A+H_B+\Omega H_c
$$
in an extended phase space which is two copies of the originial phase space. Here $H_A:=H(p,y)$ and $H_B:=H(x,q)$ correspond to two copies of the original Hamiltonian with mixed-up momentum and position variables, $H_C:=\parallel p-x \parallel_2^2/2+\parallel q-y \parallel_2^2/2$ is an imposed restraint, and $\Omega$ is a constant that controls the binding of the two copies\cite{Tao}. Then the new system is
\begin{eqnarray}\nonumber
\bordermatrix{&\cr
&\dot{p}\cr
&\dot{q}\cr
&\dot{x}\cr
&\dot{y}\cr
}
&=&\bordermatrix{&\cr
&O_n&-K_{12}(p,q)&O_n&O_n\cr
&K_{12}^{\top}(p,q)&O_n&O_n&O_n\cr
&O_n&O_n&O_n&-K_{12}(x,y)\cr
&O_n&O_n&K_{12}^{\top}(x,y)&O_n\cr
}\bordermatrix{&\cr
&\partial_p\bar{H}(p,q,x,y)\cr
&\partial_q\bar{H}(p,q,x,y)\cr
&\partial_x\bar{H}(p,q,x,y)\cr
&\partial_y\bar{H}(p,q,x,y)\cr
}\\ \label{newsys1}
&=&\bar{K}\nabla \bar{H}.
\end{eqnarray}

The original initial value problem is
\begin{align}\label{origsys}
\left\{
\begin{array}{lll}
\dot{P}=-K_{12}(P,Q)\partial_Q H(P,Q),\quad P(0)=P_0,\\
\dot{Q}=K_{12}^{\top}(P,Q)\partial_P H(P,Q),\quad Q(0)=Q_0,
\end{array}
\right.
\end{align}
and the extended system is
\begin{align}\label{augsys}
\left\{
\begin{array}{lll}
\dot{p}=-K_{12}(p,q)\partial_q H(x,q)-K_{12}(p,q)\Omega(q-y),\quad p(0)=P(0),\\
\dot{q}=K_{12}^{\top}(p,q)\partial_p H(p,y)-K_{12}^{\top}(p,q)\Omega (p-x),\quad q(0)=Q(0),\\
\dot{x}=-K_{12}(x,y)\partial_y H(p,y)-K_{12}(x,y)\Omega(y-q),\quad x(0)=P(0),\\
\dot{y}=K_{12}^{\top}(x,y)\partial_x H(x,q)+K_{12}^{\top}(x,y)\Omega(x-p),\quad y(0)=Q(0).
\end{array}
\right.
\end{align}
It should be noted that if $p(t)=x(t)$ in the extended system (\ref{augsys}), then the extended system is just the two copies of the original system, thus we have the same exact solution $p(t)=x(t)=P(t)$, $q(t)=y(t)=Q(t)$ if the same initial values are imposed on $p,x,P$ and $q,y,Q$.

Denote the time-$\tau$ flow of $H_A$ and $H_B$ by $\phi_{H_A}^{\tau}$ and $\phi_{H_B}^{\tau}$. $H_C$ can be separated into two parts $H_C^1=\parallel p-x\parallel_2^2/2$ and $H_C^{2}=\parallel q-y \parallel_2^2/2$. We denote the exact solution of the subsystem with Hamiltonian $H_C^1$ by $\phi_{\Omega H_C^1}^{\tau}$ and that with $H_C^2$ by $\phi_{\Omega H_C^2}^{\tau}$. Then we have the following four subsystems:

subsystem I
\begin{align}
\label{subsystem1}
\left\{
\begin{array}{lll}
\dot{p}=0\\
\dot{q}=K_{12}^{\top}(p,q)\partial_p H(p,y)\\
\dot{x}=-K_{12}(x,y)\partial_y H(p,y)\\
\dot{y}=0,
\end{array}
\right.
\end{align}
subsystem II
\begin{align}
\label{subsystem2}
\left\{
\begin{array}{lll}
\dot{p}=-K_{12}(p,q)\partial_q H(x,q)\\
\dot{q}=0\\
\dot{x}=0\\
\dot{y}=K_{12}^{\top}(x,y)\partial_x H(x,q),
\end{array}
\right.
\end{align}
subsystem III
\begin{align}
\label{subsystem3}
\left\{
\begin{array}{lll}
\dot{p}=0\\
\dot{q}=K_{12}^{\top}(p,q)\Omega(p-x)\\
\dot{x}=0\\
\dot{y}=K_{12}^{\top}(x,y)\Omega(x-p),
\end{array}
\right.
\end{align}
and subsystem IV
\begin{align}
\label{subsystem4}
\left\{
\begin{array}{lll}
\dot{p}=-K_{12}(p,q)\Omega(q-y)\\
\dot{q}=0\\
\dot{x}=-K_{12}(x,y)\Omega(y-q)\\
\dot{y}=0.
\end{array}
\right.
\end{align}
The problem is that under which situation, the four subsystems can be solved explicitly and the explicit K-symplectic methods can be constructed. For the above problem, we have the following theorem.

\begin{theorem}\label{theorem1}
In the following situation, the above four subsystems (\ref{subsystem1}), (\ref{subsystem2}), (\ref{subsystem3}) and (\ref{subsystem4}) can be solved explicitly:

 If the elements $K_{12}=diag(g_1(p_1,q_1),g_2(p_2,q_2),\cdots,g_n(p_n,q_n))$, there exist continuous functions $F_i(p_i,q_i)$, $G_i(p_i,q_i)$ such that $\frac{dF_i(p_i,q_i)}{dp_i}=\frac{1}{g_i(p_i,q_i)}$ and $\frac{dG_i(p_i,q_i)}{dq_i}=\frac{1}{g_i(p_i,q_i)}$, and $p_i$ and $q_i$ can be solved in closed forms from the two functions $F_i(p_i,q_i)$, $G_i(p_i,q_i)$ for each $i=1,2,\cdots,n$\cite{Zhu2}.
\end{theorem}
\proof
For the proof we refer to the details in References\cite{Zhu2}.

If $\phi_{H_A}^{\tau},\phi_{H_B}^{\tau}, \phi_{\Omega H_C^1}^{\tau}$ and $\phi_{\Omega H_C^2}^{\tau} $ can be written in close forms, then the method $\phi_{\Omega H_C^2}^{\tau}\circ \phi_{\Omega H_C^1}^{\tau}\circ \phi_{H_B}^{\tau} \circ \phi_{H_A}^{\tau}$ is a first order explicit K-symplectic method. Higher order explicit K-symplectic methods can be constructed by composing the first order K-symplectic method. The method $\phi_{ H_A}^{\tau/2}\circ \phi_{H_B}^{\tau/2} \circ \phi_{\Omega H_C^1}^{\tau/2}\circ\phi_{\Omega H_C^2}^{\tau}\circ \phi_{\Omega H_C^1}^{\tau/2}\circ \phi_{H_B}^{\tau/2} \circ \phi_{ H_A}^{\tau/2}$ is a second order explicit K-symplectic method\cite{Strang}.

\begin{remark}
For arbitrary two dimensional nonseparable non-canonical Hamiltonian system, the four subsystems can be solved exactly. Thus the K-symplectic methods can be constructed in the extended phase space by composing the exact solution of the four subsystems.
\end{remark}

For high dimensional non-separable systems, if the matrix $K^{-1}$ has the special structure as (\ref{Kstructure}) and satisfies the requirements in the above Theorem \ref{theorem1}, then we only need to make two copies of the phase space to construct explicit K-symplectic methods. If the matrix does not have the special structure, we will show in the next section another way to construct explicit K-symplectic methods.

\section{Explicit K-symplectic methods for general nonseparable systems}
\label{algorithm2}

Consider the $d$ dimensional non-canonical Hamiltonian system (\ref{noncanHam}) with nonseparable Hamiltonian $H(P)=H(p_1,p_2,\cdots,p_d)$. The skew-symmetric matrix $K^{-1}(Z)$ has no special structure. Denote by $K^{-1}=(k_{ij})_{d\times d}$. We extend the $d$ dimensional phase space to $d^2$ dimensional phase space. Denote by $P_1=(p_{11},p_{12},\cdots,p_{1d})$, $P_2=(p_{21},p_{22},\cdots,p_{2d}), \cdots, P_d=(p_{n1},p_{n2},\cdots,p_{dd})$.

We first consider the augmented Hamiltonian with mix-up variables as follows
$$
\bar{H}(P_1,\cdots,P_d)=\sum_{i=1}^d H(p_{d+1-i,1},p_{d+2-i,2},\cdots,p_{d,i},p_{1,i+1},\cdots,p_{d-i,d})+\Omega H_c.
$$
Denote by $H_1=H(p_{d1},p_{12},p_{23},\cdots,p_{d-1,d})$, $H_2=H(p_{d-1,1},p_{d2},p_{13},\cdots,p_{d-2,d})$, $\cdots$, $H_{d-1}=H(p_{21},p_{32},\cdots,p_{d,d-1},p_{1d})$ and
 $H_d=H(p_{11},p_{22},\cdots,p_{d-1,d-1},p_{dd})$,
then $\bar{H}=H_1+H_2+\cdots+H_d+\Omega H_c$ where for any $1\le i\le d$, $H_i$ take one variable from $P_j$, $1\le j\le d$. The function $H_c$ is a constraint with
%\begin{eqnarray*}
%H_c&=&\frac{\parallel P_1-P_2 \parallel_2^2}{2}+\frac{\parallel P_1-P_3 \parallel_2^2}{2}\cdots+\frac{\parallel P_1-P_{d-1} \parallel_2^2}{2}+\frac{\parallel P_1-P_{d} \parallel_2^2}{2}\\
%& &+\frac{\parallel P_2-P_3 \parallel_2^2}{2}+\frac{\parallel P_2-P_4 \parallel_2^2}{2}\cdots+\frac{\parallel P_2-P_d \parallel_2^2}{2}\\
%& &+\cdots\cdots\cdots\\
%& &+\frac{\parallel P_{d-1}-P_d \parallel_2^2}{2}.
%\end{eqnarray*}
\begin{eqnarray*}
H_c=\sum_{j=2}^d\frac{\parallel P_1-P_j \parallel_2^2}{2}+\sum_{j=3}^d \frac{\parallel P_2-P_j \parallel_2^2}{2}+\cdots+\sum_{j=d-1}^d \frac{\parallel P_{d-2}-P_j \parallel_2^2}{2}+\frac{\parallel P_{d-1}-P_d \parallel_2^2}{2}.
\end{eqnarray*}

Then the new system is
\begin{eqnarray}\nonumber
\bordermatrix{&\cr
&\dot{P}_1\cr
&\dot{P}_2\cr
&\vdots\cr
&\dot{P}_d\cr
}
&=&\bordermatrix{&\cr
& K^{-1}(P_1)& & &  \cr
& &K^{-1}(P_2)& & \cr
& & & \ddots & \cr
& & & &K^{-1}(P_d)\cr
}\bordermatrix{&\cr
&\partial_{P_1}\bar{H}(P_1,P_2,\cdots,P_d)\cr
&\partial_{P_2}\bar{H}(P_1,P_2,\cdots,P_d)\cr
&\vdots\cr
&\partial_{P_d}\bar{H}(P_1,P_2,\cdots,P_d)\cr
}\\\label{newsys2}
&=&\bar{K}\nabla \bar{H}.
\end{eqnarray}
It can be verified that if the conditions $P_1(t)=P_2(t)=\cdots=P_d(t)$ hold, the above system (\ref{newsys2}) is the $d$ copies of the original system (\ref{noncanHam}), thus we have $P_1(t)=P_2(t)=\cdots=P_d(t)=P(t)$ if the same initial condition is imposed on $P$ and each $P_i, 1\le i\le d$.

The next step is to explicitly solve each subsystem with the Hamiltonian $H_i$ for $1\le i\le d$. Here we present the subsystem with $H_1=H(p_{d1},p_{12},p_{23},\cdots,p_{d-1,d})$
\begin{align}
\label{subsystem11}
\left\{
\begin{array}{lll}
\dot{p}_{1i}=k_{i,2}(P_1)\frac{\partial H_1}{\partial p_{12}},\quad i=1,3,4,\cdots,d\\
\dot{p}_{12}=0\\
\dot{p}_{2i}=k_{i,3}(P_2)\frac{\partial H_1}{\partial p_{23}},\quad i=1,2,4,\cdots,d\\
\dot{p}_{23}=0\\
\cdots\cdots\cdots\\
\dot{p}_{d1}=0\\
\dot{p}_{d,i}=k_{i,1}(P_{d})\frac{\partial H_1}{\partial p_{d1}},\quad i=2,3,\cdots,d.\\
\end{array}
\right.
\end{align}
It can be seen that $p_{12},p_{23},\cdots,p_{d1}$ are all constants as their derivatives are zero. As the Hamiltonian $H_1$ is just a function of constants $p_{12},p_{23}$, $\cdots,p_{d1}$, therefore its partial derivatives are all constants. Thus, whether the $d$ subsystems can be solved explicitly or not depends on $k_{ij}, i,j=1,\cdots,d$.

For the constraint function $\Omega H_c$, we separate it into $\Omega H_c=\sum_{i=d+1}^{d^2}\Omega H_{i}$ with $H_{d+i}=\sum_{j=2}^d \frac{(p_{1i}-p_{ji})^2}{2}$, $H_{2d+i}=\sum_{j=3}^d \frac{(p_{2i}-p_{ji})^2}{2}$, $\cdots$, $H_{(d-1)d+i}=\frac{(p_{d-1,i}-p_{di})^2}{2}, i=1,2,\cdots,d$. Here we only present the subsystem with Hamiltonian $\Omega H_{d+1}$
\begin{align}
\label{subsystem1n}
\left\{
\begin{array}{lll}
\dot{p}_{11}=0\\
\dot{p}_{1i}=k_{i1}(P_1)\Omega (p_{11}-p_{21}+p_{11}-p_{31}+\cdots+p_{11}-p_{d1}),\quad i=2,3,\cdots,d\\
\dot{p}_{21}=0\\
\dot{p}_{2i}=k_{i1}(P_2)\Omega (p_{21}-p_{11}),\quad i=2,3,\cdots,d\\
\cdots\cdots\cdots\\
\dot{p}_{d1}=0\\
\dot{p}_{d,i}=k_{i1}(P_{d})\Omega (p_{d1}-p_{11}),\quad i=2,3,\cdots,d.\\
\end{array}
\right.
\end{align}
We can find that $p_{11},p_{21},\cdots,p_{d1}$ are all constants.

Now we discuss how to identify the situation in which the above $d^2$ subsystems can be solved explicitly. We make an explanation under the original variable $P=(p_1,p_2,\cdots,p_d)$ and consider the matrix $K^{-1}(P)=(k_{ij}(P))_{d\times d}$.
 From the subsystem (\ref{subsystem11}) and (\ref{subsystem1n}), we observe that if we want to solve each variable in a close form, the $k_{ij}$ should only contain the $i$-th variable $p_i$ and the $j$-th variable $p_j$, i.e. $k_{ij}=k_{ij}(p_i,p_j)$\cite{Zhu3}, furthermore, the following conditions should also hold: there exist continuous functions $F_{ij}(p_i,p_j)$ and $G_{ij}(p_i,p_j)$ such that $\frac{dF_{ij}(p_i,p_j)}{dp_i}=\frac{1}{k_{ij}(p_i,p_j)}$ and $\frac{dG_{ij}(p_i,p_j)}{dp_j}=\frac{1}{k_{ij}(p_i,p_j)}$ and $p_i$ and $p_j$ can be solved in closed forms from the functions $F_{ij}(p_i,p_j)$ and $G_{ij}(p_i,p_j)$ for any $1\le i,j\le d, i\neq j$. If the above requirements are satisfied, then all the subsystems can be solved explicitly.

\begin{remark}
By extending the $d$ dimensional phase space to the $d^2$ dimensional phase space, whether the subsystems can be solved exactly or not depends totally on the elements of the skew-symmetric matrix $K^{-1}$.
\end{remark}

If all the $d^2$ subsystems can be solved explicitly, then the first order explicit K-symplectic method can be constructed by composing all the exact solutions. The higher order explicit K-symplectic methods can be constructed by composing the first order K-symplectic method. Given a first order K-symplectic method $\Phi_{\tau}$ where $\tau$ represents the time stepsize, then we can compose it to obtain a higher order K-symplectic method
$$
\Psi_{\tau}\equiv \Phi_{\alpha_s\tau}\circ \Phi_{\beta_s \tau}^*\circ \cdots \circ \Phi_{\beta_2 \tau}^*\circ \Phi_{\alpha_1 \tau}\circ \Phi_{\beta_1 \tau}^*
$$
with the coefficients $\alpha_{i}=\beta_{s+1-i}, 1\le i\le s$. Here $\Phi_{\tau}^*$ is the adjoint method of $\Phi_{\tau}$.

In many real problems, the high dimensional problem can be decoupled into several lower dimensional problems, thus one can solve the lower dimensional problems. For example, the equation for the distribution function $f$ of the Vlasov-Maxwell equations can be rewritten as a $6n$ dimensional Hamiltonian system by expressing $f$ as a sum of Dirac masses\cite{He2}, and the $6n$ dimensional system can be decoupled into $n$ charged particle systems of $6$ dimension.

\section{Numerical experiments}

\subsection{Numerical methods}

All experiments were performed in MATLAB R2016b on a Windows 10 (64 bit) PC with the configuration: Intel(R) Core(TM) i7-8550U CPU 1.80 GHz and 16 GB RAM.  To demonstrate the superiority in the structure preservation of the K-symplectic methods, they are compared with higher order Runge-Kutta methods. Denote $\Phi_{\tau}$ by the first order K-symplectic method composed by the exact solution of all subsystems.

2ndKsym: the second order K-symplectic method, which is the composition of and its adjoint method
$$
\Psi_{\tau}^2\equiv \Phi_{\tau/2}^*\circ \Phi_{\tau/2}.
$$

4thKsym: the fourth order K-symplectic method which is composed by
$$
\Psi_{\tau}^4\equiv \Phi_{\alpha_5\tau}\circ \Phi_{\beta_5 \tau}^*\circ \cdots \circ \Phi_{\beta_2 \tau}^* \circ \Phi_{\alpha_1 \tau}\circ \Phi_{\beta_1 \tau}^*.
$$
The values of the parameters $\alpha_1,\beta_1,\cdots,\alpha_5,\beta_5$ are given in \cite{McLachlan}.

3rdRK: the third order explicit Runge-Kutta method (Heun method)\cite{HNL}. This method is compared with the 2ndKsym method.

5thRK: the explicit Runge-Kutta method with effective order 5 (Butcher method)\cite{Butcher2,HNL}. This method is compared with the 4thKsym method.

As the parameter $\Omega$ is used to control the several copies of the phase space, it should not be very small. Therefore we set $\Omega=20$.

\subsection{The first numerical demonstration}

To illustrate the behaviors of the K-symplectic methods, we perform the numerical simulation for a four dimensional non-canonical system. The skew-symmetric matrix is
$$K^{-1}(P)=
\begin{pmatrix}
0	  &	 0&	  \frac{1}{sec^2(xz)} & 0 \\	
0 & 0  &	 0 &	 \frac{u^2}{2sin(y)}\\
-\frac{1}{sec^2(xz)} &0 & 0 &  0\\
0&-\frac{u^2}{2sin(y)} &0& 0
\end{pmatrix},\quad P=(x,y,z,u)^{\top}$$
and the Hamiltonian is $H=(x^2+y^2+z^2)^{5/2}+yu$. As the Hamiltonian is not separable, thus we extend the phase space. It can be easily seen that the matrix $K^{-1}$ satisfies the requirements in Theorem \ref{theorem1}, thus we extend the four dimensional phase space to eight dimensional phase space. The augmented Hamiltonian is
\begin{eqnarray*}
& &\bar{H}(p_1,p_2,p_3,p_4,q_1,q_2,q_3,q_4)=H(p_1,p_2,q_3,q_4)+H(q_1,q_2,p_3,p_4)\\
& &+\Omega\Big(\frac{(p_1-q_1)^2}{2}+\frac{(p_2-q_2)^2}{2}+\frac{(p_3-q_3)^2}{2}+\frac{(p_4-q_4)^2}{2}\Big).
\end{eqnarray*}
Then the extended system can be separated into six subsystems with $H_1:=H(p_1,p_2,q_3,q_4)$, $H_2:=H(q_1,q_2,p_3,p_4)$ and $H_i:=\Omega\frac{(p_{i-2}-q_{i-2})^2}{2}, i=3,4,5,6.$
Here we only present the exact solution of the subsystem with $H_1$
\begin{align}
\label{exactsolu1}
\left\{
\begin{array}{lll}
p_i=p_{i0},\quad i=1,2\\
p_3=\frac{1}{p_{10}}\arctan\Big(\frac{\partial H_1}{\partial p_{10}}p_{10}t+tan(p_{10}p_{30})\Big)\\
p_4=\frac{1}{\frac{1}{2sin(p_{20})}\frac{\partial H_1}{\partial p_{20}}t+\frac{1}{p_{40}}}\\
q_1=\frac{1}{q_{30}}\arctan\Big(\frac{\partial H_1}{\partial q_{30}}q_{30}t+tan(q_{10}q_{30})\Big)\\
q_2=\arccos(\frac{-q_{40}^2}{2}\frac{\partial H_1}{\partial q_{40}}t+\cos(q_{20}))\\
q_i=q_{i0},\quad i=3,4.\\
\end{array}
\right.
\end{align}
As the Hamiltonian $H_1$ is the function of constants $p_1,p_2,q_1,q_2$, thus the partial derivative of $H_1$ with respect to each argument is also a constant. The other subsystems can also be solved explicitly, thus the explicit K-symplectic method can be constructed.

\begin{figure}[p]
\centering
\subfigure[ ]{
\includegraphics[scale=0.42]{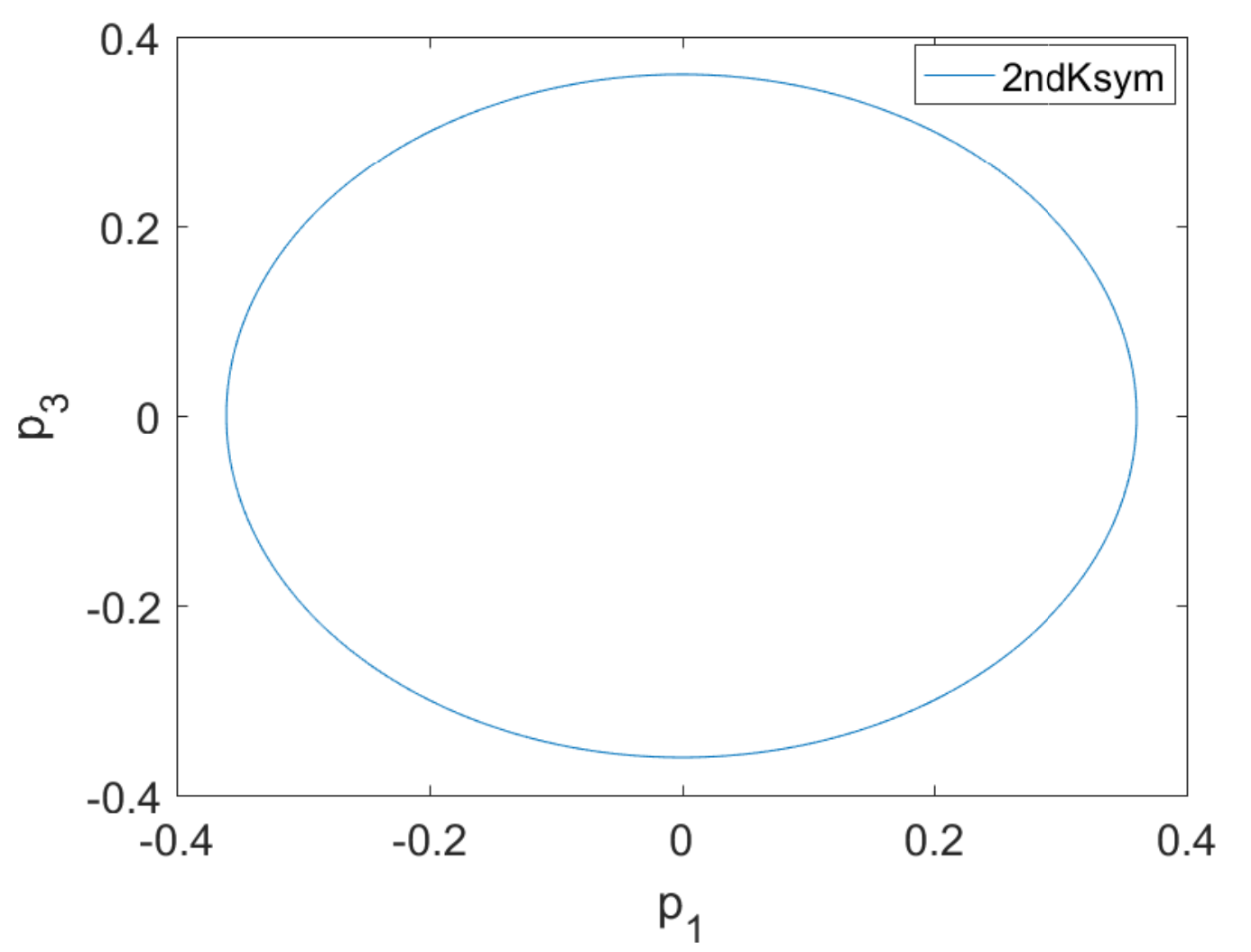}}
\subfigure[ ]{
\includegraphics[scale=0.42]{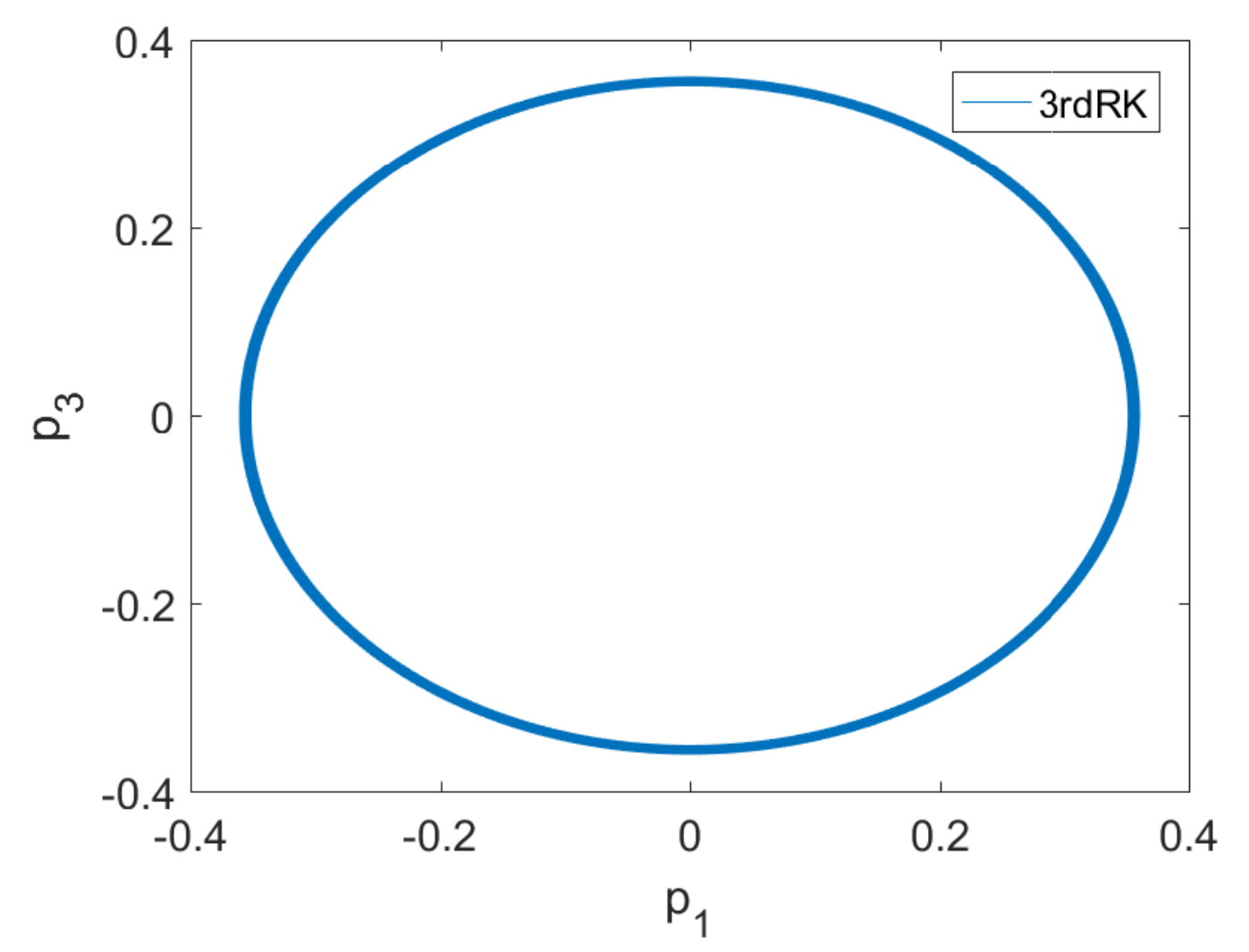}}
\subfigure[ ]{
\includegraphics[scale=0.42]{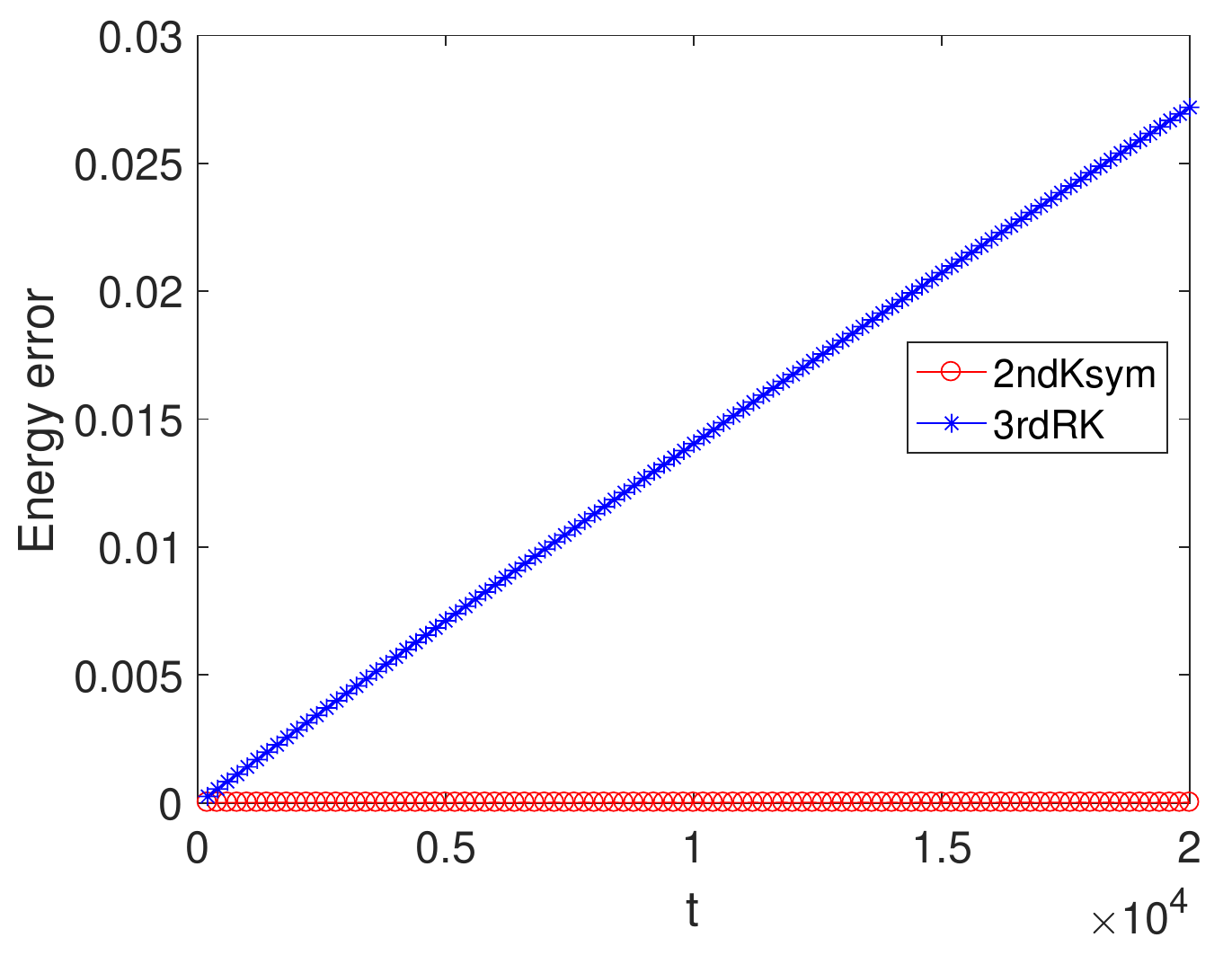}}
\subfigure[ ]{
\includegraphics[scale=0.42]{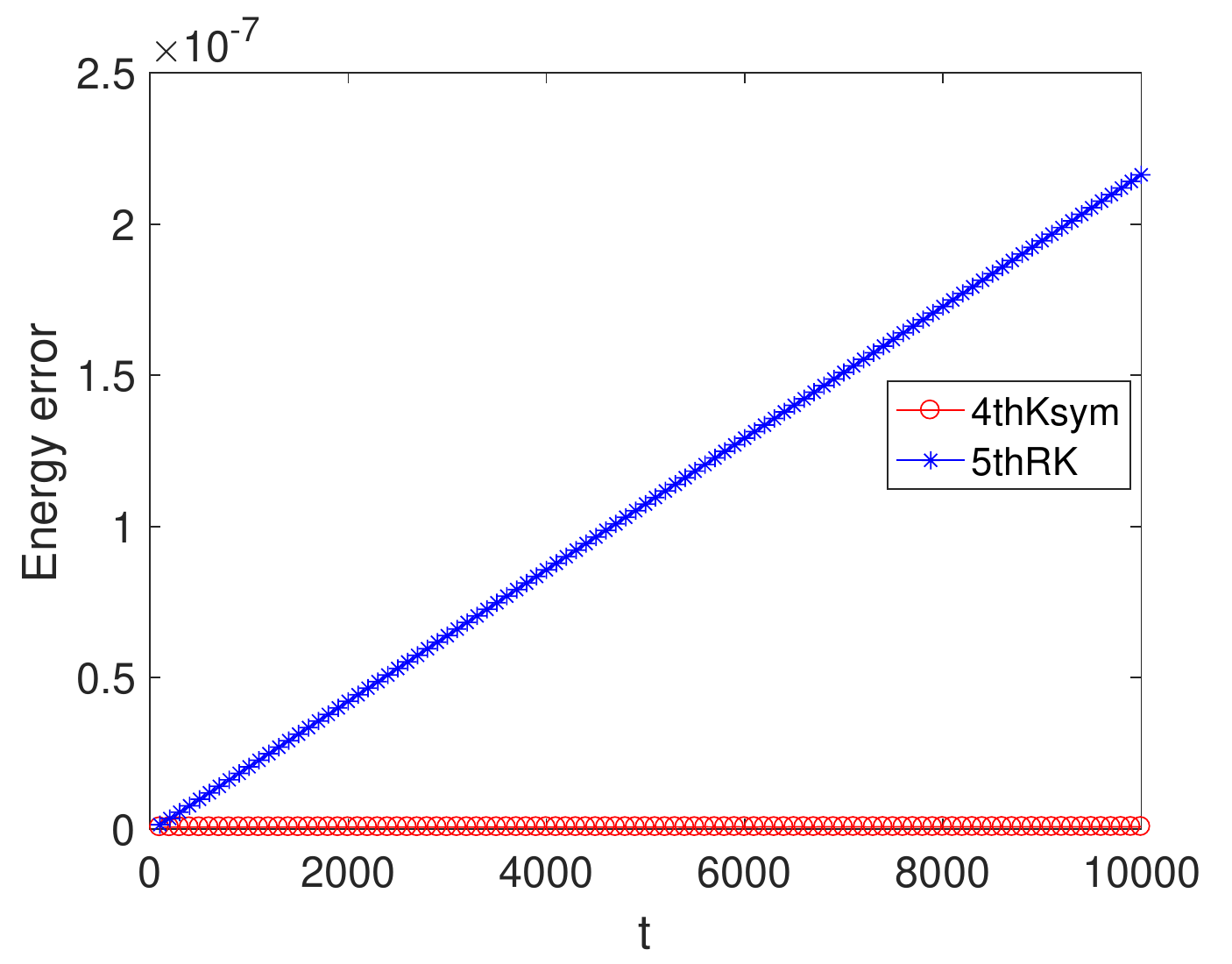}}
\caption{The phase orbit and the energy error obtained by the 2ndKsym method, the 3rdRK method, the 4thKsym method and the 5thRK method. Subfigure (a) and (b) display the orbit projected to $p_1-p_3$ plane obtained by the 2ndKsym method and the 3rdRK method with the stepsize $\tau=0.01$, $\Omega=20$ over the interval $[0,20000]$.  Subfigure (c) and (d) display the relative energy error of the augmented Hamiltonian $\bar{H}$ for the four methods with the stepsize $\tau=0.01$ and $\Omega=20$. The energy error is represented by $(\bar{H}-\bar{H}_0)/\bar{H}_0$. }
\label{fig:phase1}
\end{figure}

\begin{figure}[p]
\centering
\subfigure[ ]{
\includegraphics[scale=0.42]{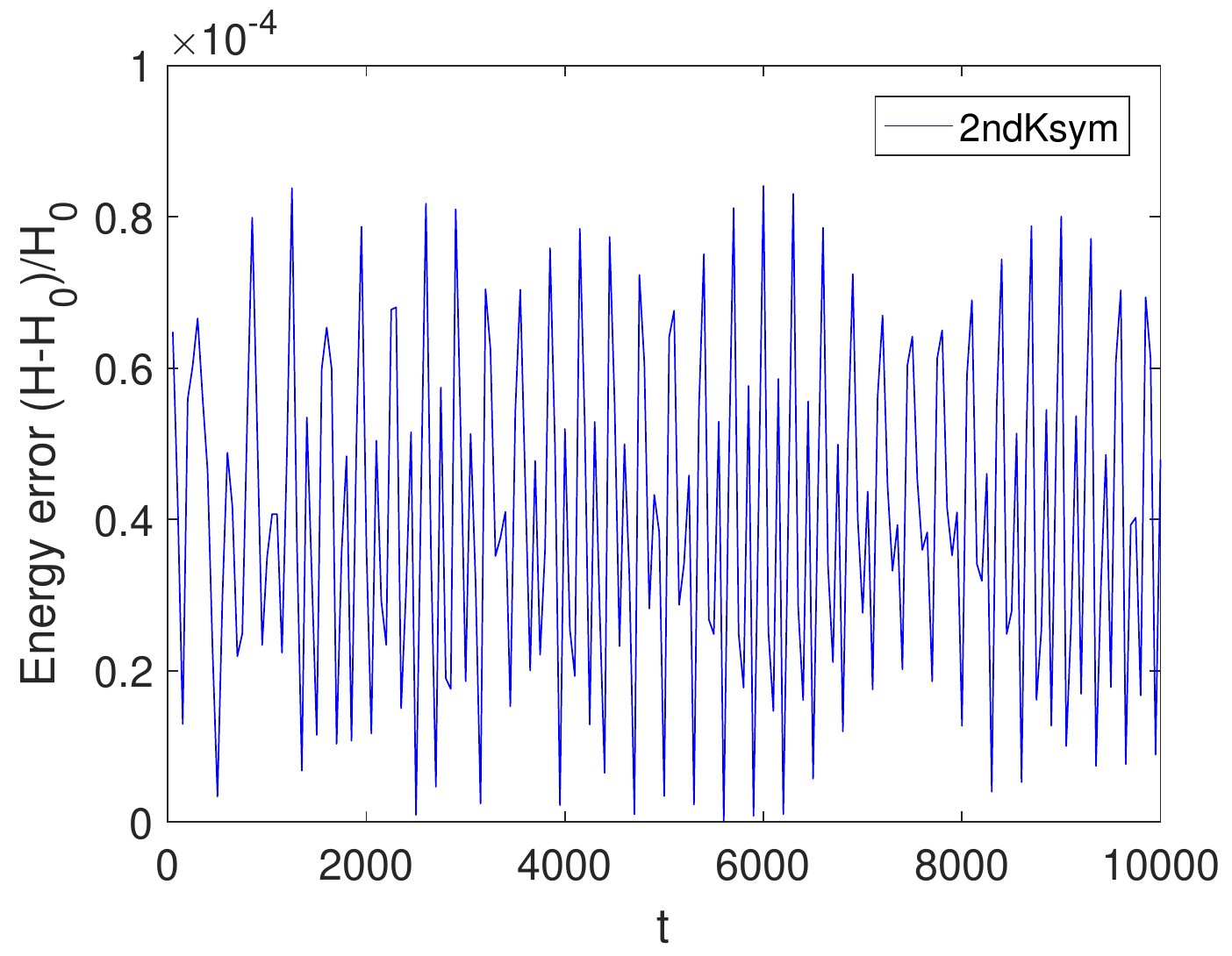}}
\subfigure[ ]{
\includegraphics[scale=0.42]{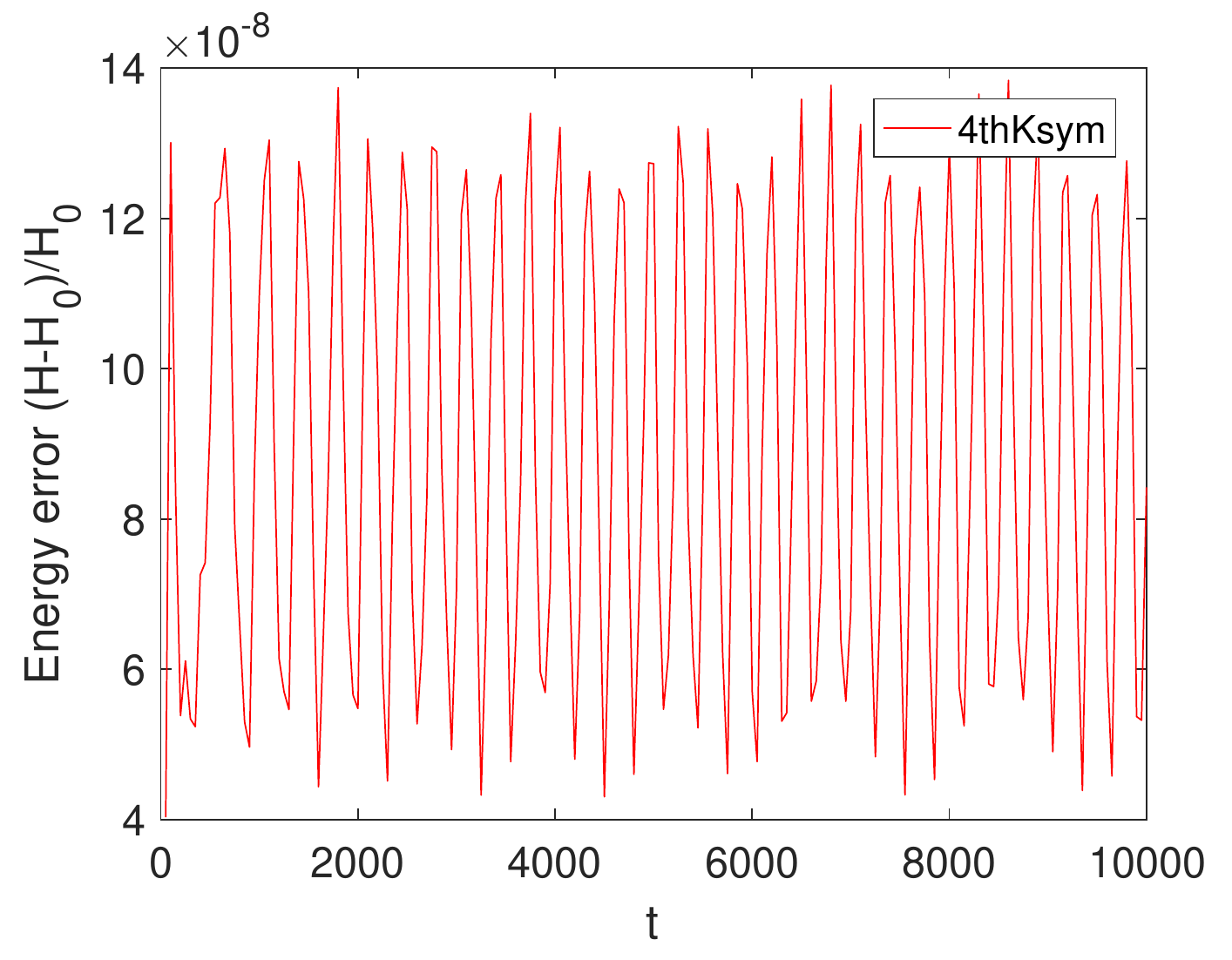}}
\caption{The relative energy error of the original Hamiltonian $H(p_1,p_2,p_3,p_4)$ obtained by the 2ndKsym method and the 4thKsym method. The relative energy error is represented by $(H-H_0)/H_0$. The stepsize is chosen as $\tau=0.01$ and the parameter is $\Omega=20$.}
\label{fig:energy1}
\end{figure}

\begin{figure}[p]
\centering
\subfigure[ ]{
\includegraphics[scale=0.42]{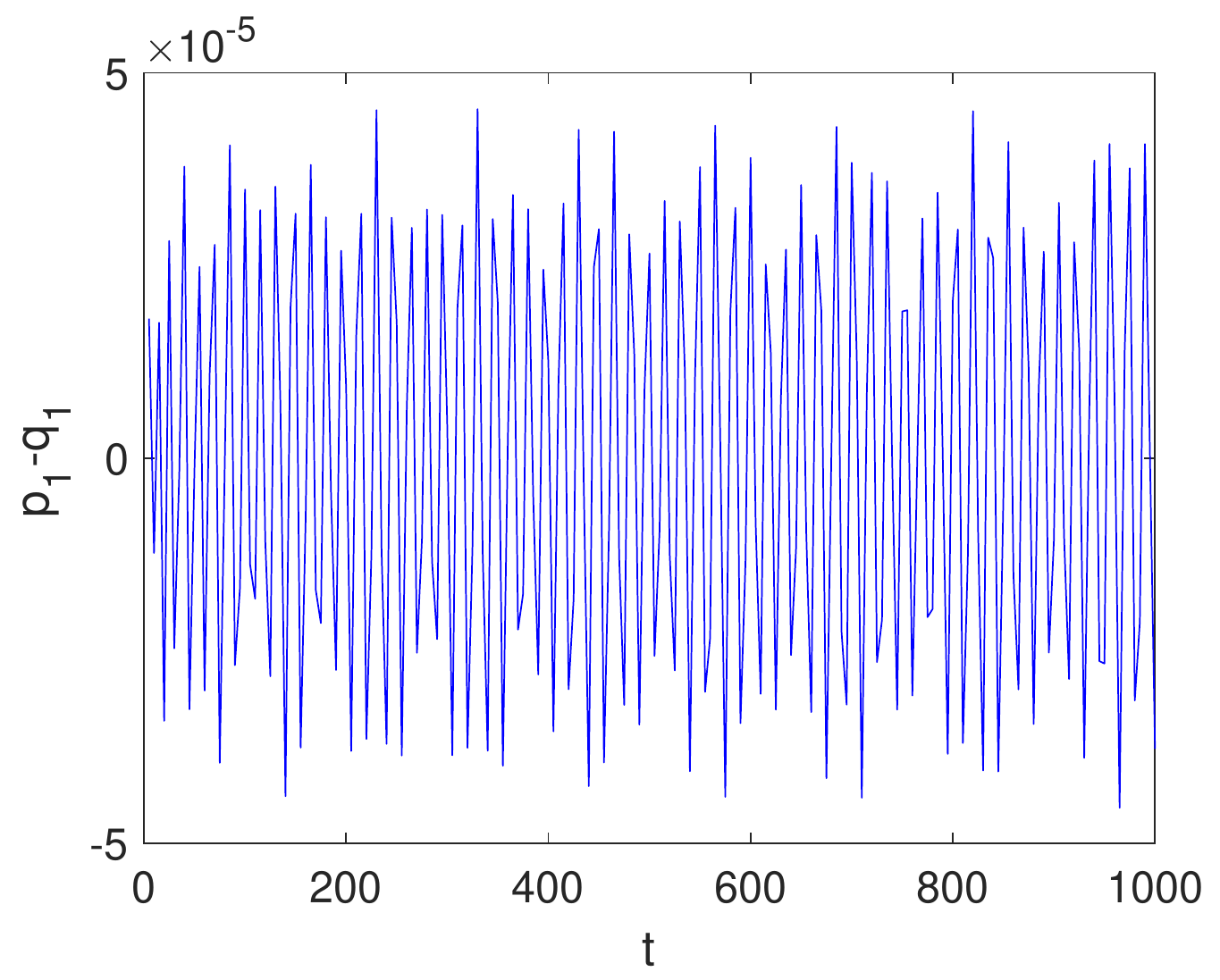}}
\subfigure[ ]{
\includegraphics[scale=0.42]{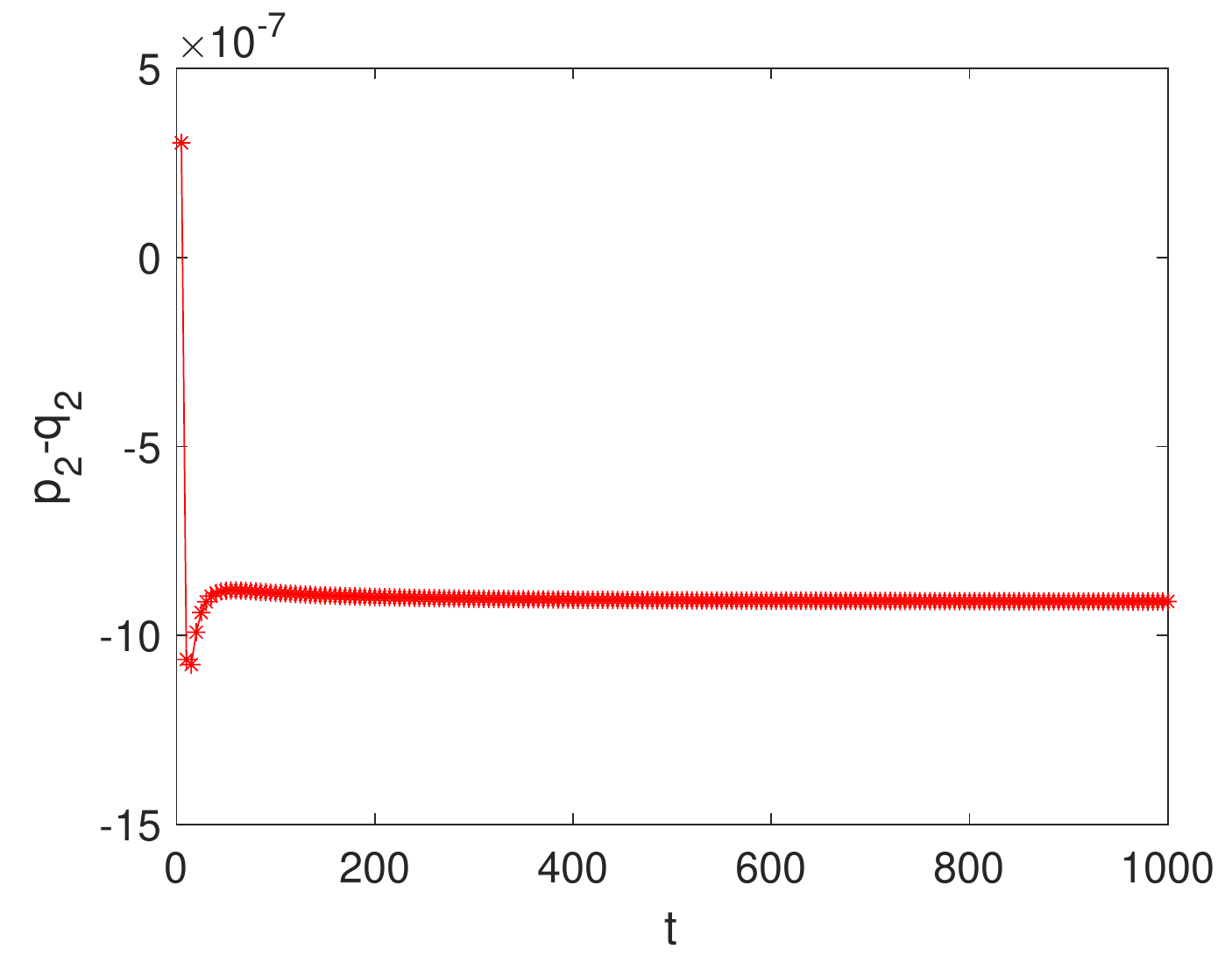}}
\subfigure[ ]{
\includegraphics[scale=0.42]{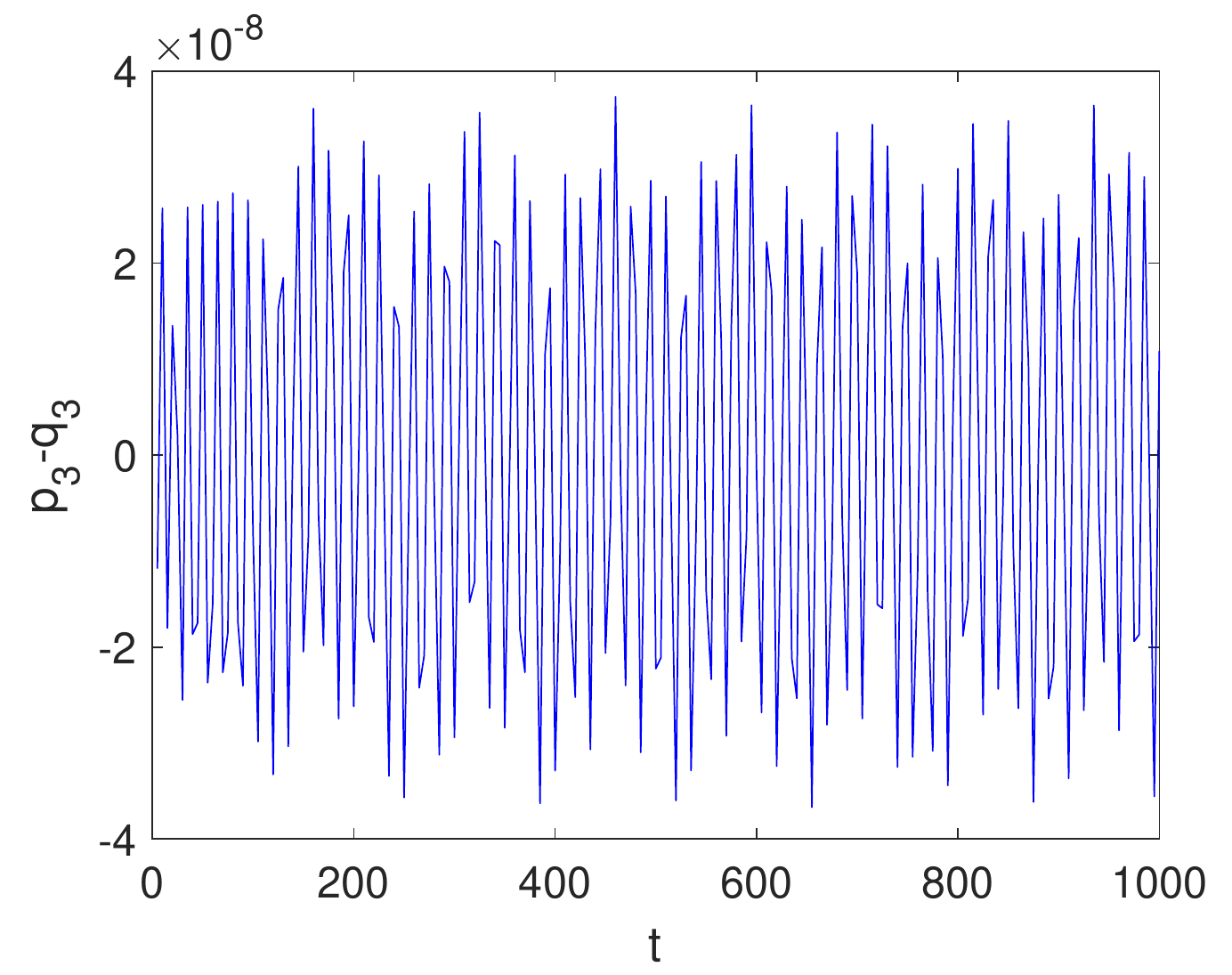}}
\subfigure[ ]{
\includegraphics[scale=0.42]{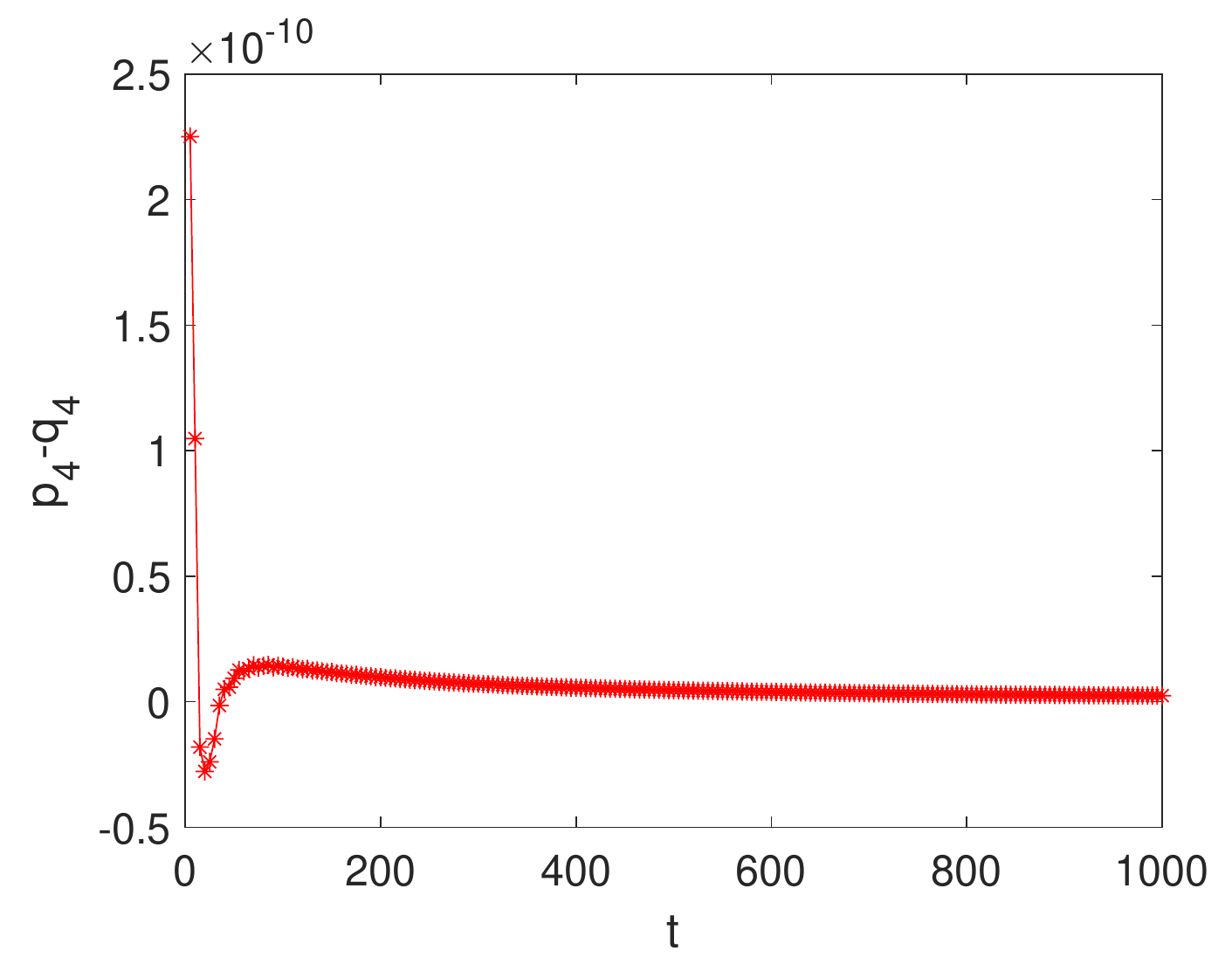}}
\caption{The difference between the two copies of the four variables. Subfigure (a), (b), (c) and (d) display the difference between the two copies for the four variables respectively.  Subfigure (a), (b) are the differences in the first two variables obtained by the 2ndKsym method. Subfigure (c), (d) are the differences in the last two variables obtained by the 4thKsym method. Here the stepsize is $\tau=0.01$ and $\Omega=20$, the final time is $T=1000$.}
\label{fig:pqerror1}
\end{figure}

The initial condition is $x_0=0.2,y_0=0.4,z_0=0.3, u_0=0.5$. Then the initial condition for the extended system is $p_{10}=q_{10}=x_0,p_{20}=q_{20}=y_0$, $p_{30}=q_{30}=z_0, p_{40}=q_{40}=u_0$. The numerical results for the four numerical methods are displayed in Figure \ref{fig:phase1}-\ref{fig:pqerror1}. The orbits in $p_1-p_3$ plane obtained by the 2ndKsym method and 3rdRK method are displayed in Figure \ref{fig:phase1}. We can see that the orbit obtained by the second order K-symplectic method is more accurate than the third order Runge-Kutta method. The evolutions of the energy $\bar{H}$ using different methods are also shown in Figure \ref{fig:phase1}. The energy errors of K-symplectic methods oscillate with small amplitudes while the energy errors of the Runge-Kutta methods increase without bound along time as can be seen from Figure \ref{fig:phase1}.
 The relative energy errors of the original $H$ for the two K-symplectic methods are shown in Figure \ref{fig:energy1}. It can be seen that the energy errors of the two methods can be bounded at a small interval. The differences between the two copies of the original variables are displayed in Fig. \ref{fig:pqerror1}. The difference in each variable is bounded at a very small number as can be seen from Figure \ref{fig:pqerror1}. The CPU times of the four methods are displayed in Table \ref{table1}. As can be seen that the computational cost of the second order explicit K-symplectic method is smaller than that of the third explicit Runge-Kutta method. The CPU time of the fourth order explicit K-symplectic method is longer than that of the fifth order explicit Runge-Kutta method, but the difference is not so big.

\begin{table}[htbp]
\begin{small}
\caption{The CPU times of the four methods. The stepsize is $\tau=0.01$ and the time interval is $[0,1000]$.}
\begin{center}
\begin{tabular}{|r|c|c|c|}
\hline
2ndKsym & 3rdRK & 4thKsym & 5thRK\\
\hline
 0.5916 & 0.7531  &  4.5615 & 1.7261 \\
\hline
\end{tabular}
\end{center}
\label{table1}
\end{small}
\end{table}

\subsection{The second numerical demonstration for Ablowitz-Ladik model}

The Ablowitz-Ladik model is one of the space discretization of the nonlinear Schr\"{o}dinger equation with the form of \cite{GNT}
\begin{equation}\label{almodel}
i\dot{W}_k+\frac{1}{h^2}(W_{k+1}-2W_k+W_{k-1})+|W_k|^2(W_{k+1}+W_{k-1})=0.
\end{equation}
where the dot represents the derivative with respect to the time $t$ and $h$ is the space stepsize. Under the periodic condition $W_{k+N}=W_k (h=\frac{1}{N})$, it is an integral system. By separating the variable $W_k$ into real and imaginary parts, $W_k=u_k+iv_k$, then we obtain a non-canonical Hamiltonian system with
$$
\bordermatrix{&\cr
&\dot{u}\cr
&\dot{v}\cr
}
=\bordermatrix{&\cr
&O &-D\cr
&D&O\cr
}\bordermatrix{&\cr
&\partial_u H(u,v)\cr
&\partial_v H(u,v)\cr
}=K^{-1}(u,v)\nabla H(u,v)
$$
where $u=(u_1,\cdots,u_N), v=(v_1,\cdots,v_n)$, $D=diag(d_1,d_2,\cdots,d_N)$ is a diagonal matrix with the elements
$$
d_k(u,v)=1+h^2(u_k^2+v_k^2),\quad k=1,2\cdots,N
$$
and the Hamiltonian is
$$
H(u,v)=\frac{1}{h^2}\sum_{i=1}^N (u_iu_{i-1}+v_iv_{i-1})-\frac{1}{h^4}\sum_{i=1}^N \ln(1+h^2(u_i^2+v_i^2)).
$$
It can be easily seen that the matrix $K^{-1}$ satisfies the requirements in
Theorem \ref{theorem1}. Thus we extend the phase space and make the two copies of $(u,v)$ and separate the Hamiltonian into two parts.

We take $N=4$ and denote the two copies of $(u,v)$ by $p=(p_1,p_2,\cdots,p_8)$ and $q=(q_1,q_2,\cdots,q_8)$. Thus we consider the augmented Hamiltonian
$$
\bar{H}:=H_1+H_2+\Omega H_c
$$
where $H_1:=H(p_1,\cdots,p_4,q_5,\cdots,q_8)$, $H_2:=H(q_1,\cdots,q_4,p_5,\cdots,p_8)$ and
$$
H_c=\frac{\parallel P-Q \parallel_2^2}{2}=\sum_{i=1}^8 \frac{(p_i-q_i)^2}{2}.
$$
Then the augmented Hamiltonian is separated into 10 subsystems, with $H_1$, $H_2$ and $H_{j+2}=\frac{(p_i-q_i)^2}{2}, j=1,2,\cdots,8$. Here we only present the exact solution of the first subsystem with the Hamiltonian $H_1$
\begin{align}
\label{almexactsolu1}
\left\{
\begin{array}{lll}
p_i=p_{i0},\quad i=1,2,3,4\\
p_j=\frac{\sqrt{1+h^2p_{j-4,0}^2}}{h}\tan\Big(t h\frac{\partial H_1}{\partial p_{j-4,0}}\sqrt{1+h^2p_{j-4,0}^2}+\arctan\Big(\frac{hp_{j0}}{\sqrt{1+h^2p_{j-4,0}^2}}\Big)\Big),j=5,6,7,8\\
q_j=\frac{\sqrt{1+h^2q_{j+4,0}^2}}{h}\tan\Big(-t h\frac{\partial H_1}{\partial q_{j+4,0}}\sqrt{1+h^2q_{j+4,0}^2}+\arctan\Big(\frac{hq_{j0}}{\sqrt{1+h^2q_{j+4,0}^2}}\Big)\Big), \quad i=1,2,3,4\\
q_i=q_{i0},\quad i=5,6,7,8.\\
\end{array}
\right.
\end{align}
where $\tau$ is the time stepsize and $h$ is the space stepsize.
As the Hamiltonian $H_1$ is the function of constants $p_{i0},q_{i+4,0}, i=1,2,3,4$, therefore their partial derivatives with respect to these arguments are all constants. The other subsystems can also be solved explicitly, therefore the explicit K-symplectic methods can be constructed.

\begin{figure}[p]
\centering
\subfigure[ ]{
\includegraphics[scale=0.42]{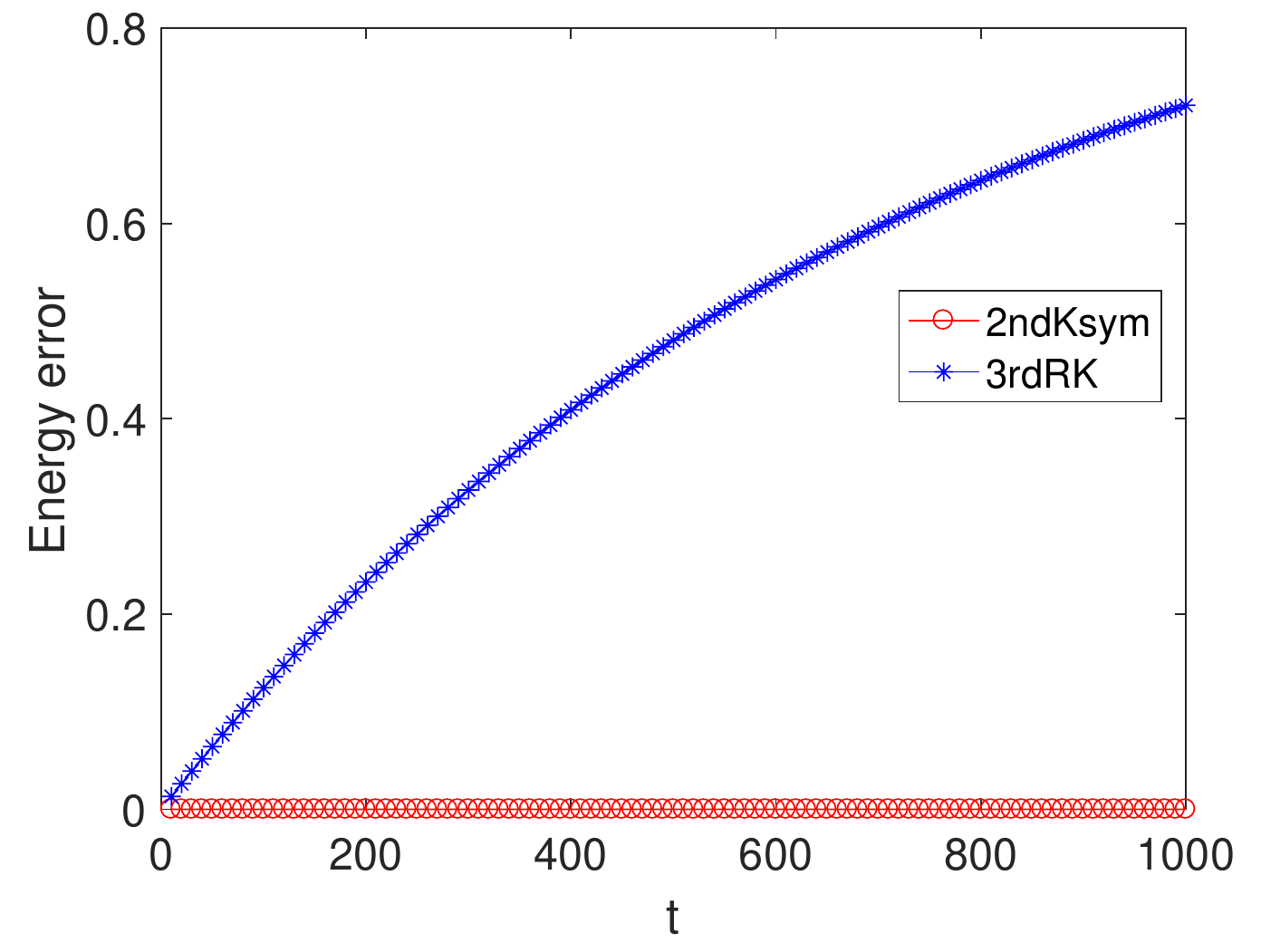}}
\subfigure[ ]{
\includegraphics[scale=0.42]{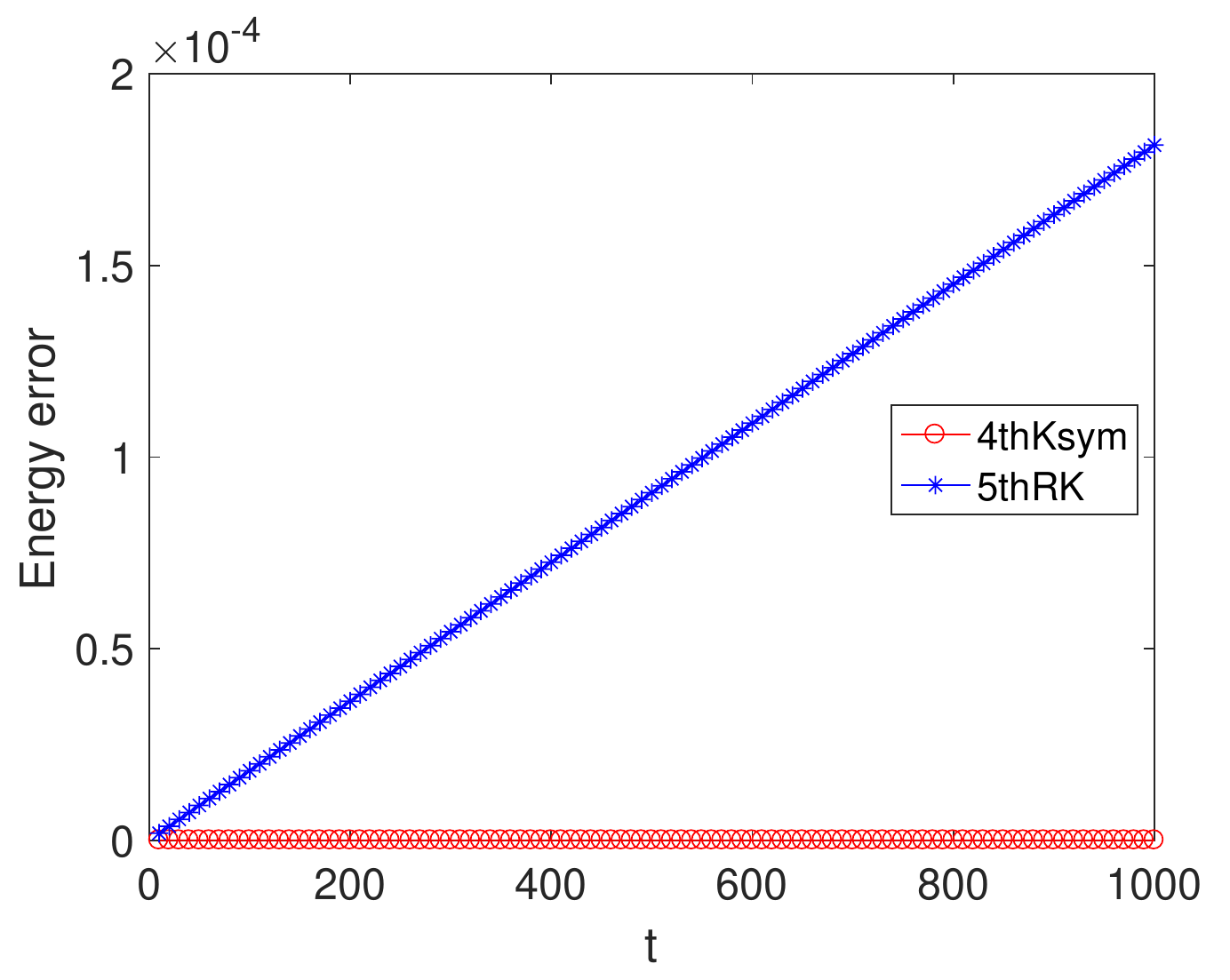}}
\subfigure[ ]{
\includegraphics[scale=0.42]{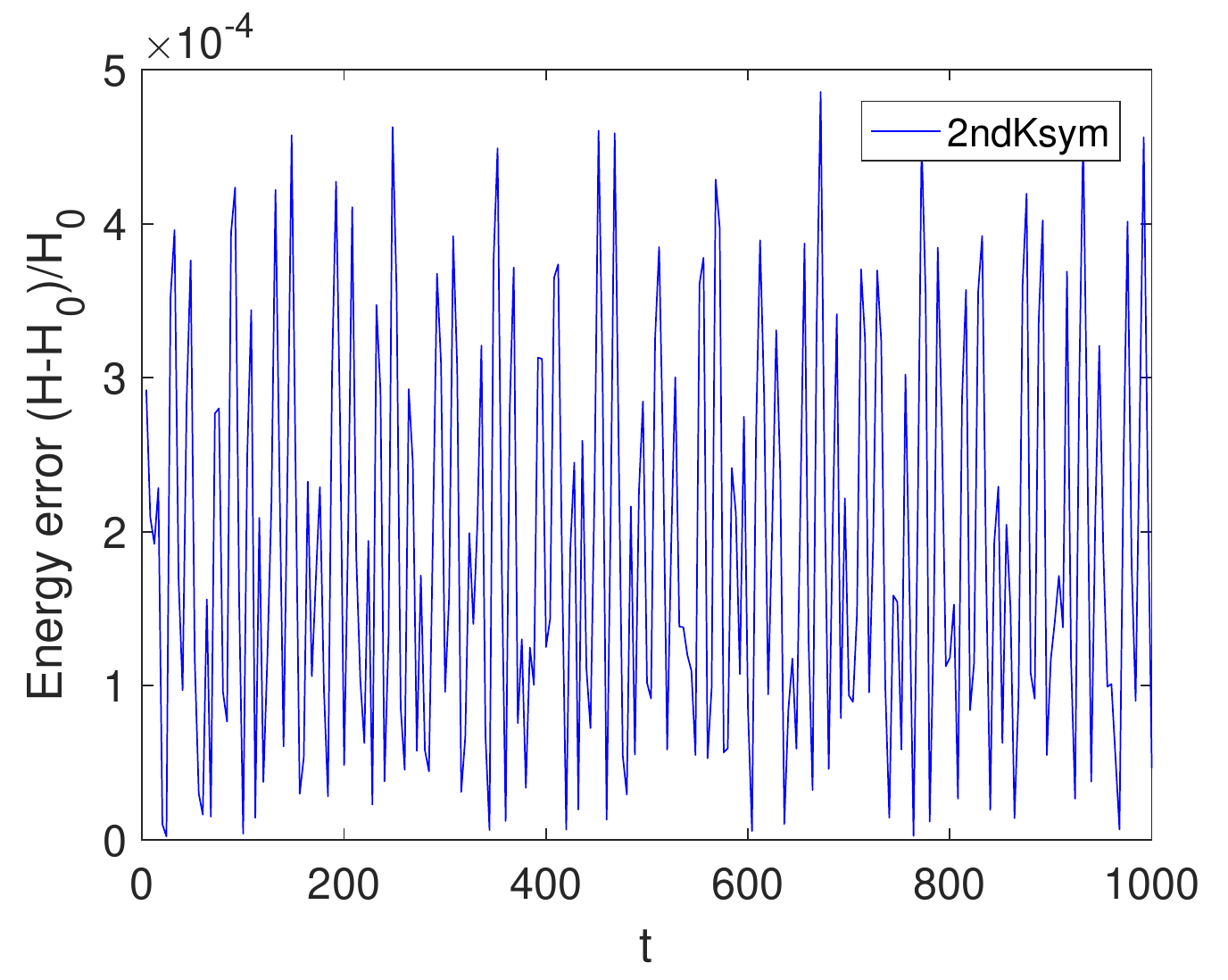}}
\subfigure[ ]{
\includegraphics[scale=0.42]{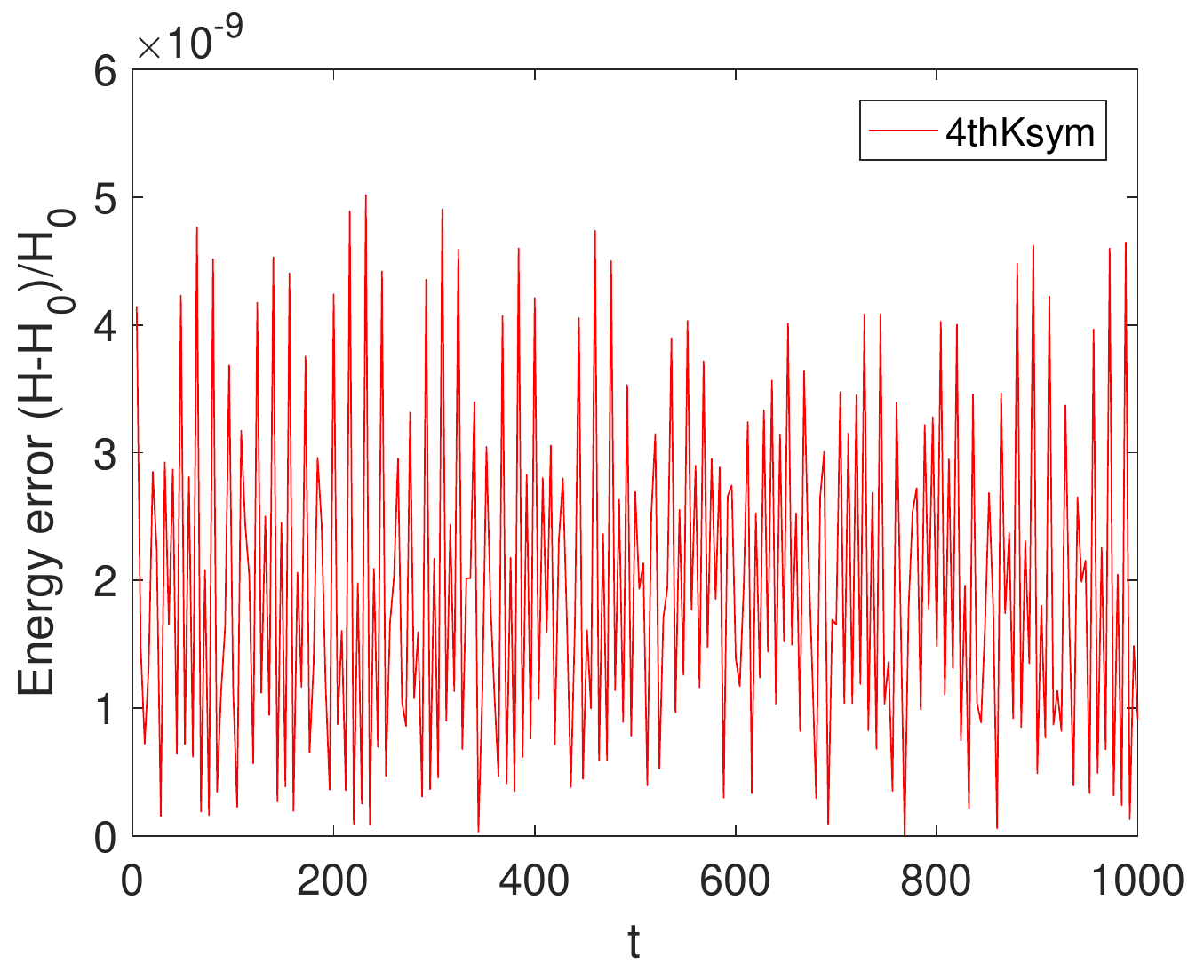}}
\caption{The energy evolutions of the 2ndKsym method, the 3rdRK method, the 4thKsym method and the 5thRK method.  Subfigure (a) and (b) display the relative energy error of the augmented Hamiltonian $\bar{H}$ for the four methods with $\Omega=20$ and the time stepsize $\tau=0.001$. In the first two subfigures the energy error is represented by $(\bar{H}-\bar{H}_0)/\bar{H}_0$. Subfigure (a) and (b) display the relative energy error of the original Hamiltonian $H$ with $\Omega=20$ and $\tau=0.001$. In the last two subfigures the energy error is represented by $(H-H_0)/H_0$.}
\label{fig:energy3}
\end{figure}

\begin{figure}[p]
\centering
\subfigure[ ]{
\includegraphics[scale=0.42]{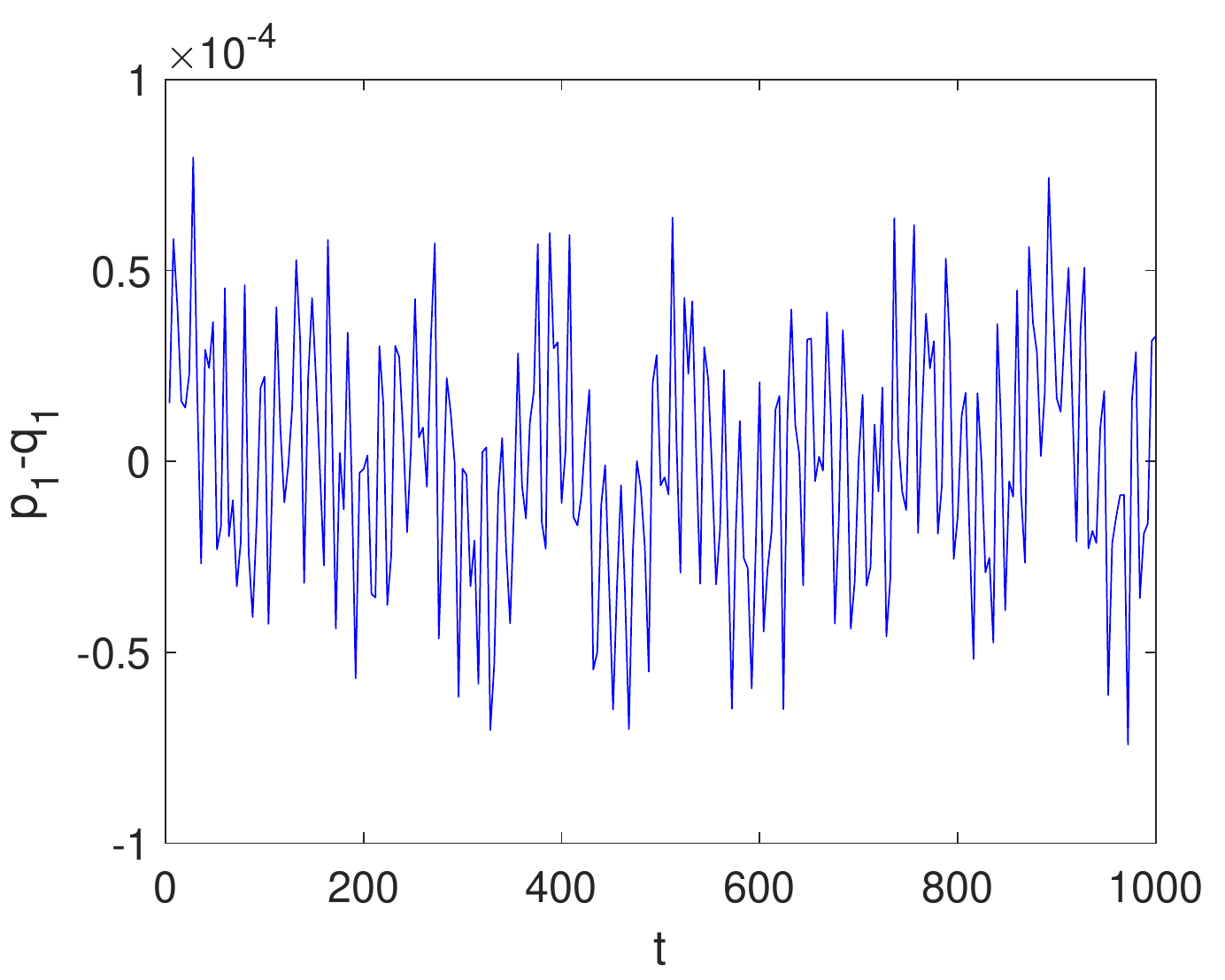}}
\subfigure[ ]{
\includegraphics[scale=0.42]{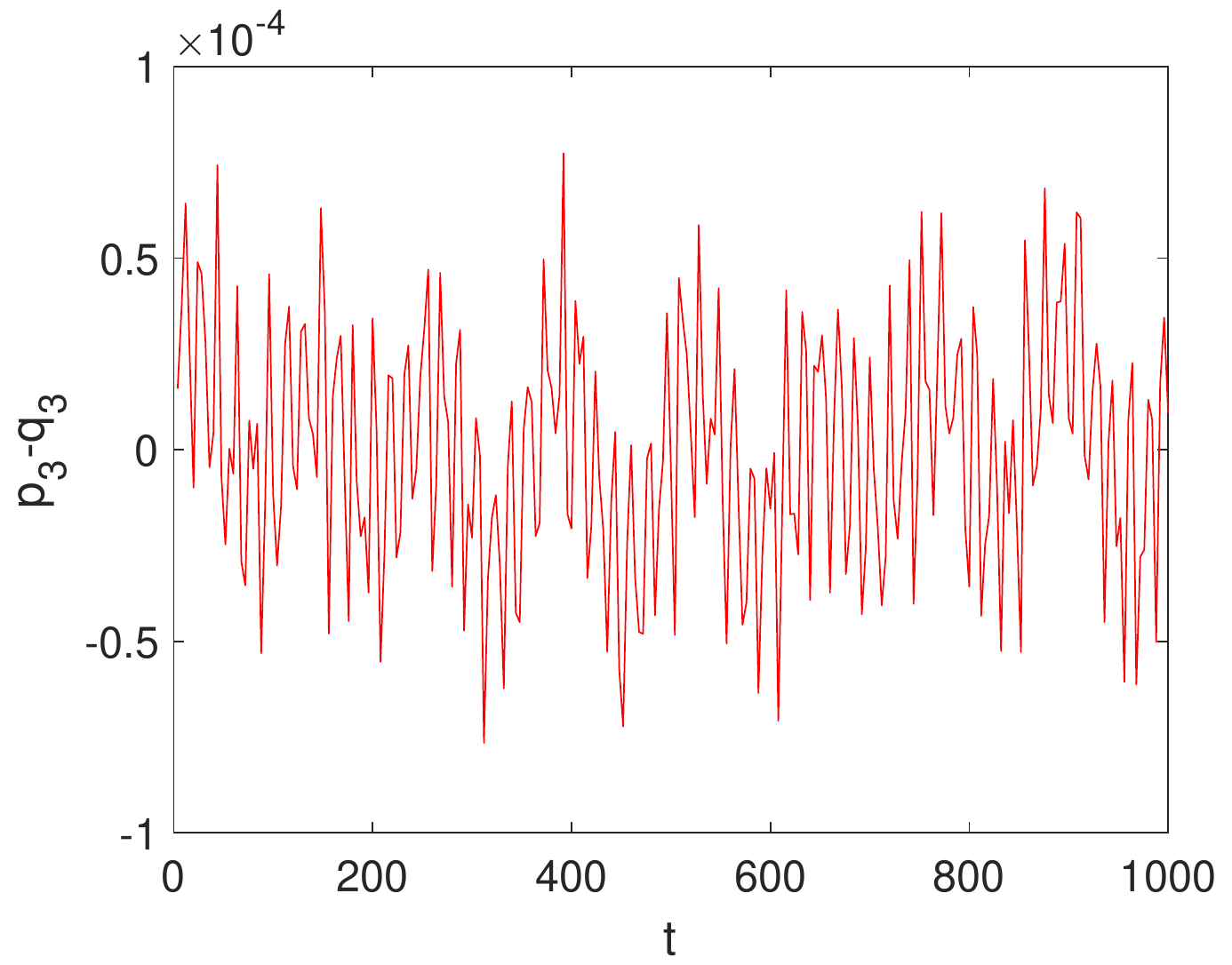}}
\subfigure[ ]{
\includegraphics[scale=0.42]{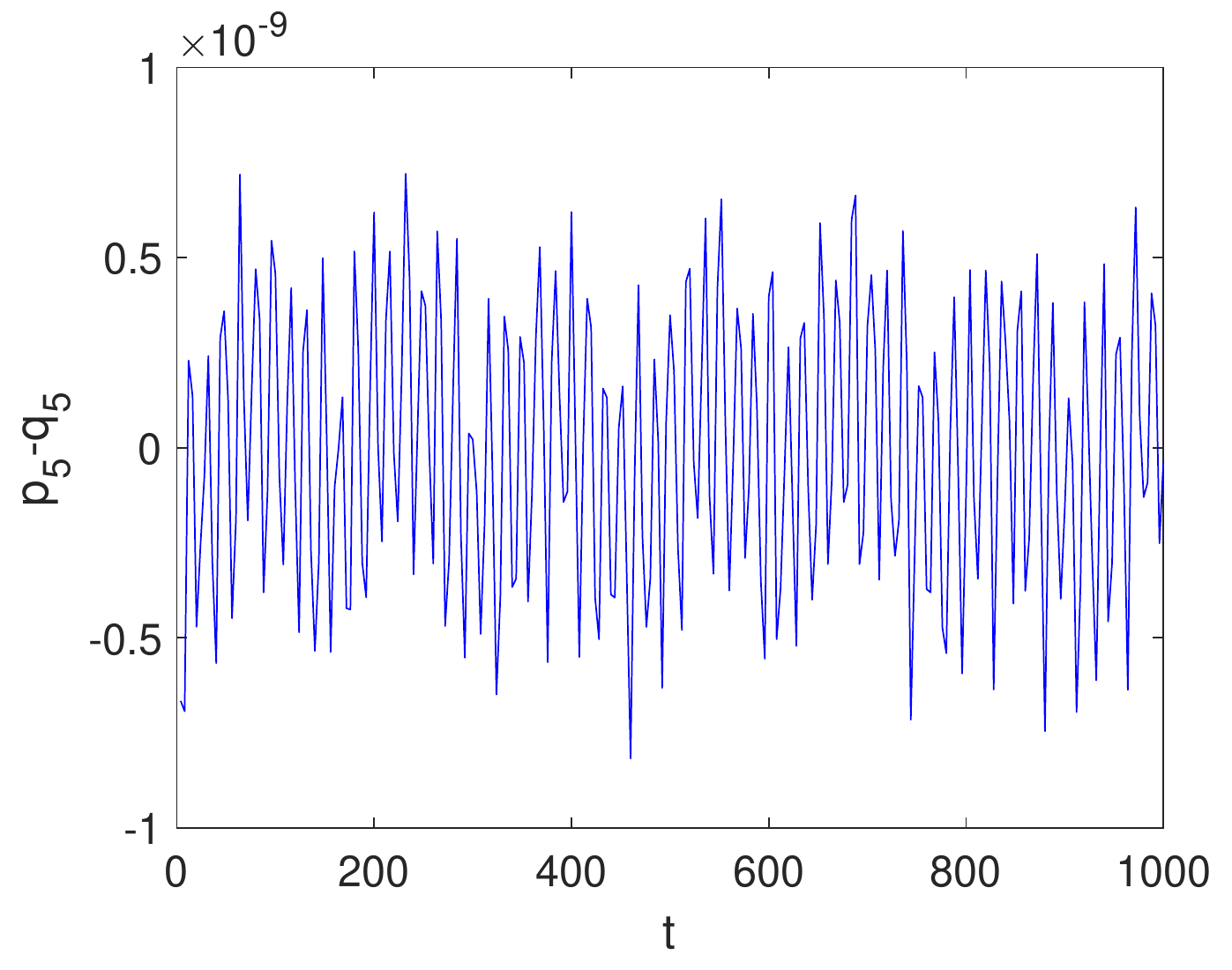}}
\subfigure[ ]{
\includegraphics[scale=0.42]{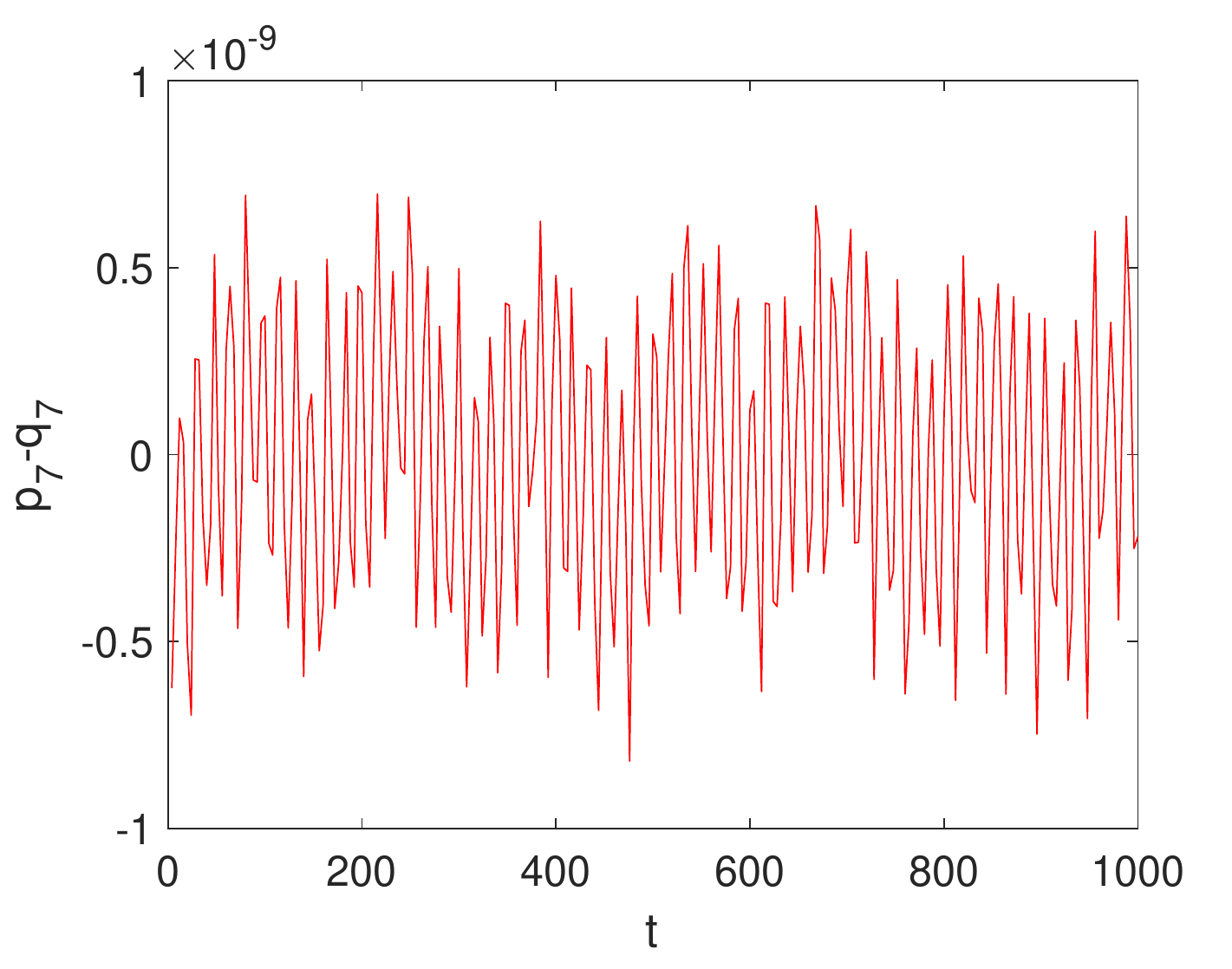}}
\caption{The difference between the two copies for four variables $u_1$, $u_3$, $v_1$ and $v_3$. Subfigure (a) and (b) displays the differences between the two copies for the variables $u_1$ and $u_3$ obtained by the 2ndKsym method. Subfigure (c) and (d) displays the differences between the two copies for the variables $v_1$ and $v_3$ obtained by the 4thKsym method. Here the time stepsize is $\tau=0.001$ and $\Omega=20$, the final time is $T=1000$.}
\label{fig:pqerror3}
\end{figure}

The initial condition is $u_{10}=0.2,u_{20}=0.4,u_{30}=0.3,u_{40}=0.5$, $v_{10}=0.3,v_{20}=0.2$, $v_{30}=0.3,v_{40}=0.2$. The two copies of the variable $(u,v)$ have the same initial condition. The energy evolutions of the augmented Hamiltonian $\bar{H}$ obtained by the 2ndKsym, 4thKsym, 3rdRK and 5thRK methods are displayed in Figure \ref{fig:energy3}. The explicit K-symplectic methods have shown their significant advantages in energy conservation over long-term simulation compared with the higher order Runge-Kutta method. The energy errors of the original Hamiltonian $H(p_1,\cdots,p_8)$ in the first copy of variables obtained by the two explicit K-symplectic methods are also shown in Figure \ref{fig:energy3}.  As can be seen from Figure \ref{fig:energy3} that the energy errors oscillate with a small amplitude. The two copies of the variable $(u,v)$ are also compared. We show in Figure \ref{fig:pqerror3} the differences between the two copies of $u_1,u_3,v_1,v_3$. The errors in the two copies of $u_1,u_3,v_1,v_3$ obtained by the two explicit K-symplectic methods are all bounded by small numbers.

\subsection{The third numerical demonstration for gyrocenter system}

Gyrocenter system\cite{Littlejohn,Qin,Zhang,Zhu3} plays a fundamental role in plasma physics research. It is a non-canonical Hamiltonian system
\begin{equation}\label{origsys2}
\dot{Z}=K^{-1}(Z)\nabla H(Z),\quad Z=(x,y,z,u)^{\top},
\end{equation}
with the Hamiltonian $H=\frac{1}{2}u^{2}+\mu |B(X)|+\varphi({X})$ and
$$K^{-1}(Z)=\dfrac{1}{a_{12}b_3-a_{13}b_2+a_{23}b_1}
\begin{pmatrix}
0 & -b_3 & b_2 & a_{23} \\ b_3 & 0 & -b_1 & -a_{13} \\
-b_2 & b_1 & 0 & a_{12} \\ -a_{23} & a_{13} & -a_{12} & 0
\end{pmatrix}.$$
Given a vector potential $A(X)=(f,g,h)^\top$, the magnetic field is $B(X)=\nabla\times A(X)$, $b(X)=(b_1,b_2,b_3)^\top=\dfrac{B(X)}{|B(X)|}$ and the elements in $K^{-1}(Z)$ are
$$a_{12}=g_x-f_y+u(\dfrac{\partial b_2}{\partial x}-\dfrac{\partial b_1}{\partial y}),$$ $$a_{13}=h_x-f_z+u(\dfrac{\partial b_3}{\partial x}-\dfrac{\partial b_1}{\partial z}),$$ $$a_{23}=h_y-g_z+u(\dfrac{\partial b_3}{\partial y}-\dfrac{\partial b_2}{\partial z}).$$

We set the magnetic field is $B(X)=(0,0,sec^2(xy))$. There exist continuous functions $f$ and $g$ satisfying that $g_x-f_y=sec^2(xy)$. The scalar potential is set to be $\varphi(X)=10^{-2}\sqrt{x^2+y^2+z^2}$. The Hamiltonian is $H=sec^2(xy)+10^{-2}\sqrt{x^2+y^2+z^2}+u^2/2$. As the Hamiltonian is not separable and $K^{-1}$ does not have the special structure in Section \ref{Kspecial}, therefore we make four copies of the variable $(x,y,z,u)$.

We extend the phase space $(x,y,z,u)$ to sixteen dimensional phase space.
Denote the four copies of $(x,y,z,u)$ by $P=(p_1,p_2,p_3,p_4)$, $Q=(q_1,q_2,q_3,q_4)$, $R=(r_1,r_2,r_3,r_4)$ and $W=(w_1,w_2,w_3,w_4)$, then we consider an augmented Hamiltonian
$$
\bar{H}(P,Q,R,W):=H_1+H_2+H_3+H_4+\Omega H_c,
$$
where $H_1:=H(p_1,q_2,r_3,w_4)$, $H_2:=H(w_1,p_2,q_3,r_4)$, $H_3:=H(r_1,w_2,p_3,q_4)$, $H_4:=H(q_1,r_2,w_3,p_4)$ and
\begin{eqnarray*}
H_c:&=&\frac{\parallel P-Q \parallel_2^2}{2}+\frac{\parallel P-R \parallel_2^2}{2}+\frac{\parallel P-W \parallel_2^2}{2}+\frac{\parallel Q-R \parallel_2^2}{2}\\
& &+\frac{\parallel Q-W \parallel_2^2}{2}+\frac{\parallel R-W \parallel_2^2}{2}.
\end{eqnarray*}
As the augmented Hamiltonian $\bar{H}$ is separable, we separate the original system into several subsystems. The first four subsystems are with the Hamiltonian $H_i, i=1,2,3,4$. Denote by $K^{-1}=(k_{ij})_{4\times 4}$.
%Here we only present the first subsystem
%\begin{align}
%\label{foursubsystem1}
%\left\{
%\begin{array}{lll}
%\dot{p}_1=0\\
%\dot{p}_i=k_{i1}(p_1,p_2,p_3,p_4)\frac{\partial H_1(p_1,q_2,r_3,w4)}{\partial p_1}, \quad i=2,3,4\\
%\dot{q}_i=k_{i2}(q_1,q_2,q_3,q_4)\frac{\partial H_1(p_1,q_2,r_3,w4)}{\partial q_2}, \quad i=1,3,4\\
%\dot{q}_2=0\\
%\dot{r}_i=k_{i3}(r_1,r_2,r_3,r_4)\frac{\partial H_1(p_1,q_2,r_3,w4)}{\partial r_3}, \quad i=1,2,4\\
%\dot{r}_3=0\\
%\dot{w}_i=k_{i4}(w_1,w_2,w_3,w_4)\frac{\partial H_1(p_1,q_2,r_3,w4)}{\partial w_4}, \quad i=1,2,3\\
%\dot{w}_4=0.\\
%\end{array}
%\right.
%\end{align}
%As the time derivatives of $p_1, q_2, r_3$ and $w_4$ are all zero, then we know that $\frac{\partial H_1}{\partial p_1}$, $\frac{\partial H_1}{\partial q_2}$, $\frac{\partial H_1}{\partial r_3}$ and $\frac{\partial H_1}{\partial w_4}$ are all constants.
The exact solution to the first subsystem can be solved explicitly
\begin{align}
\label{gyroexactsolu1}
\left\{
\begin{array}{lll}
p_i=p_{i0},\quad i=1,3,4\\
p_2=\frac{1}{p_{10}}\arctan\Big(\frac{\partial H_1}{\partial p_{10}}p_{10}t+tan(p_{10}p_{20})\Big)\\
q_1=\frac{1}{q_{20}}\arctan\Big(-\frac{\partial H_1}{\partial q_{20}}q_{10}t+tan(q_{10}q_{20})\Big)\\
q_i=q_{i0}, \quad i=2,3,4\\
r_i=r_{i0},\quad i=1,2,3\\
r_4=r_{40}-\frac{\partial H_1}{\partial r_{30}}t\\
w_i=w_{i0},\quad i=1,2,4\\
w_3=w_{30}+\frac{\partial H_1}{\partial w_{40}}t.\\
\end{array}
\right.
\end{align}
where $p_{i0},q_{i0},r_{i0},w_{i0}$ represent the initial values of $p_i,q_i,r_i,w_i$ for $i=1,2,3,4$. As the Hamiltonian $H_1$ is the function of the constants $p_1,q_2,r_3,w_4$, therefore the partial derivative of $H$ with respect to each argument is also a constant.

Now we are concerned with the further partition of $\Omega H_c$. We separate the $\Omega H_c$ into 12 parts, i.e. $\Omega H_c=\sum_{i=5}^{16}H_i$ with the following functions
$$
H_{4+i}=\Omega \Big[ \frac{(p_i-q_i)^2}{2}+ \frac{(p_i-r_i)^2}{2} +\frac{(p_i-w_i)^2}{2}\Big], \quad i=1,2,3,4,
$$
$$
H_{8+j}=\Omega \Big[ \frac{(q_j-r_j)^2}{2}+ \frac{(q_j-w_j)^2}{2}\Big],\quad j=1,2,3,4,
$$
and
$$
H_{12+k}=\Omega \Big[ \frac{(r_k-w_k)^2}{2}\Big], \quad k=1,2,3,4.
$$
%Here we only present the subsystem with $H_5$
%\begin{align}
%\label{foursubsystem5}
%\left\{
%\begin{array}{lll}
%\dot{p}_1=0\\
%\dot{p}_i=k_{i1}(p_1,p_2,p_3,p_4)\Omega((p_1-q_1)+(p_1-r_1)+(p_1-w_1)), \quad i=2,3,4\\
%\dot{q}_1=0\\
%\dot{q}_i=k_{i1}(q_1,q_2,q_3,q_4)\Omega(q_1-p_1), \quad i=2,3,4\\
%\dot{r}_1=0\\
%\dot{r}_i=k_{i1}(r_1,r_2,r_3,r_4)\Omega(r_1-p_1), \quad i=2,3,4\\
%\dot{w}_1=0\\
%\dot{w}_i=k_{i1}(w_1,w_2,w_3,w_4)\Omega(w_1-p_1), \quad i=2,3,4.\\
%\end{array}
%\right.
%\end{align}
%As the time derivatives of $p_1, q_1, r_1$ and $w_1$ are all zero, then we know that $\Omega((p_1-q_1)+(p_1-r_1)+(p_1-w_1))$, $\Omega(q_1-p_1)$, $\Omega(r_1-p_1)$ and $\Omega(w_1-p_1)$ are all constants. We can also solve it explicitly
Here we present the exact solution to the subsystem with Hamiltonian $H_5$
\begin{align}
\label{exactsolu5}
\left\{
\begin{array}{lll}
p_i=p_{i0},\quad i=1,3,4\\
p_2=\frac{1}{p_{10}}\arctan\Big(\Omega (p_{10}-q_{10}+p_{10}-r_{10}+p_{10}-w_{10})p_{10}t+tan(p_{10}p_{20})\Big)\\
q_i=q_{i0},\quad i=1,3,4\\
q_2=\frac{1}{q_{10}}\arctan\Big(\Omega (q_{10}-p_{10})q_{10}t+tan(q_{10}q_{20})\Big)\\
r_i=r_{i0},\quad i=1,3,4\\
r_2=\frac{1}{r_{10}}\arctan\Big(\Omega (r_{10}-p_{10})r_{10}t+tan(r_{10}r_{20})\Big)\\
w_i=w_{i0},\quad i=1,3,4\\
w_2=\frac{1}{w_{10}}\arctan\Big(\Omega (w_{10}-p_{10})w_{10}t+tan(w_{10}w_{20})\Big).\\
\end{array}
\right.
\end{align}
where $p_{i0},q_{i0},r_{i0},w_{i0}$ represent the initial values of $p_i,q_i,r_i,w_i$ for $i=1,2,3,4$.
As the Hamiltonian $H_5$ is the function of the constants $p_1,q_1,r_1,w_1$, therefore the partial derivative of $H$ with respect to each argument is also a constant.

As all the subsystems can be solved explicitly, thus explicit K-symplectic methods can be constructed by composing the exact solution of all the subsystems.

\begin{figure}[p]
\centering
\subfigure[ ]{
\includegraphics[scale=0.42]{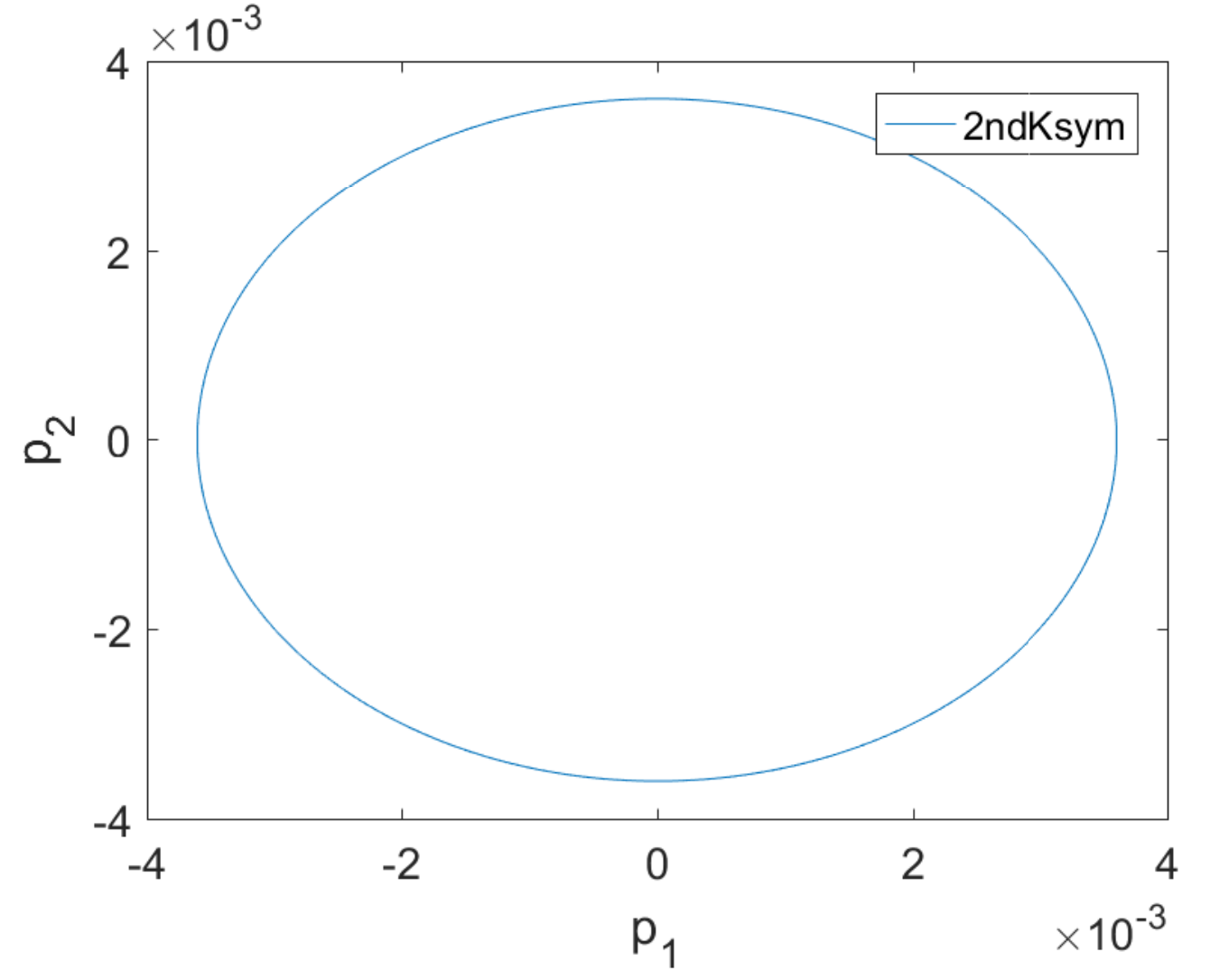}}
\subfigure[ ]{
\includegraphics[scale=0.42]{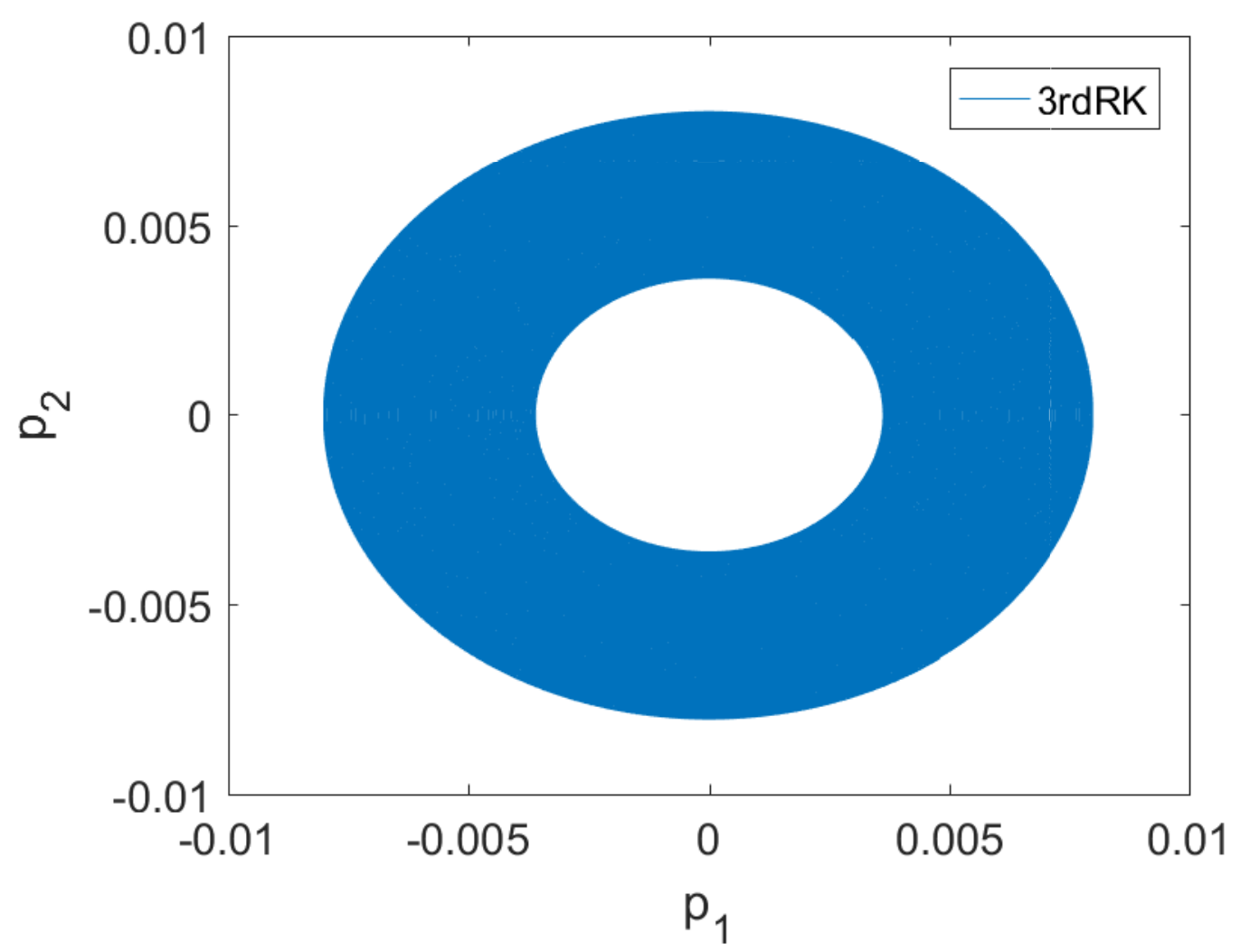}}
\subfigure[ ]{
\includegraphics[scale=0.42]{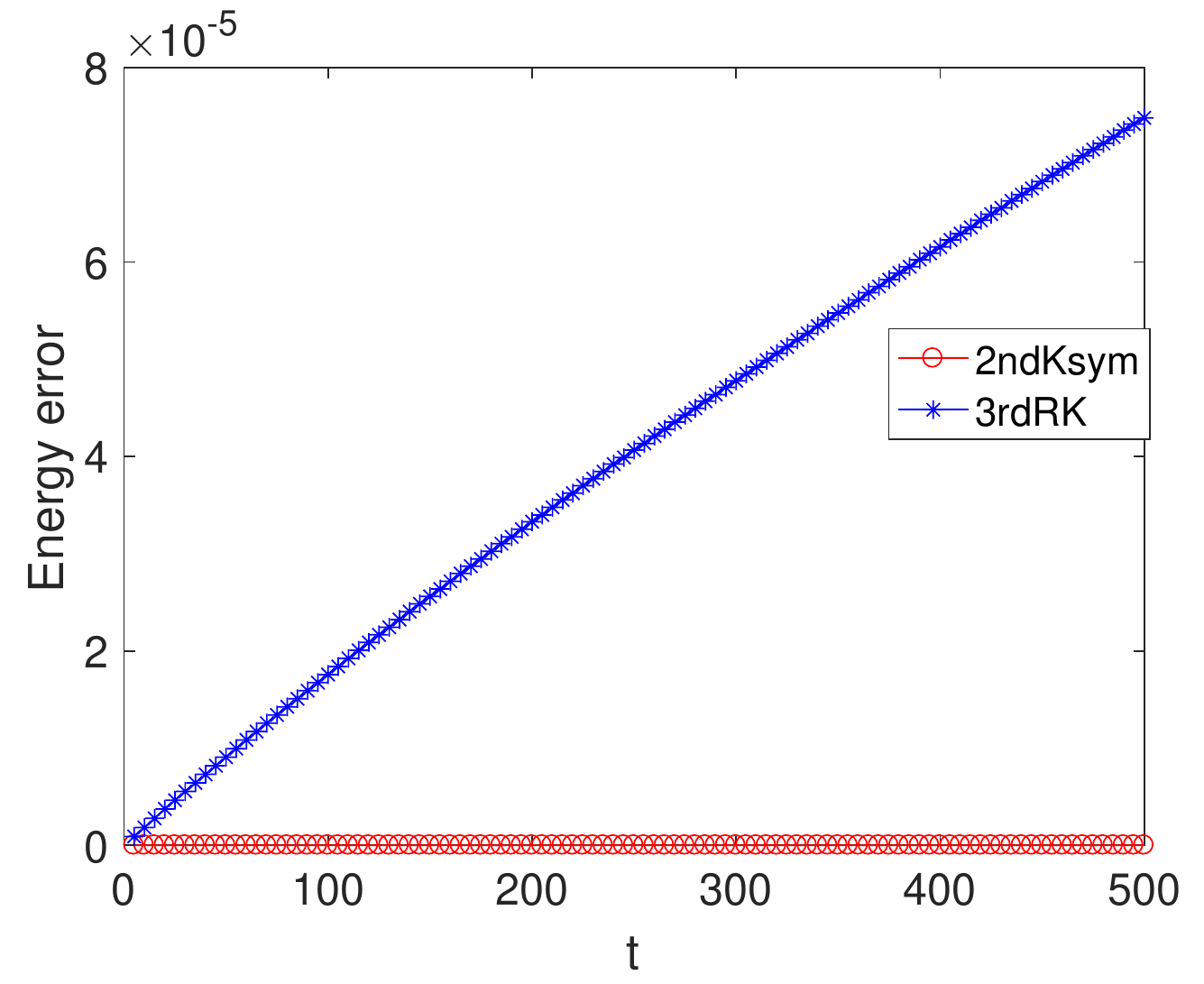}}
\subfigure[ ]{
\includegraphics[scale=0.42]{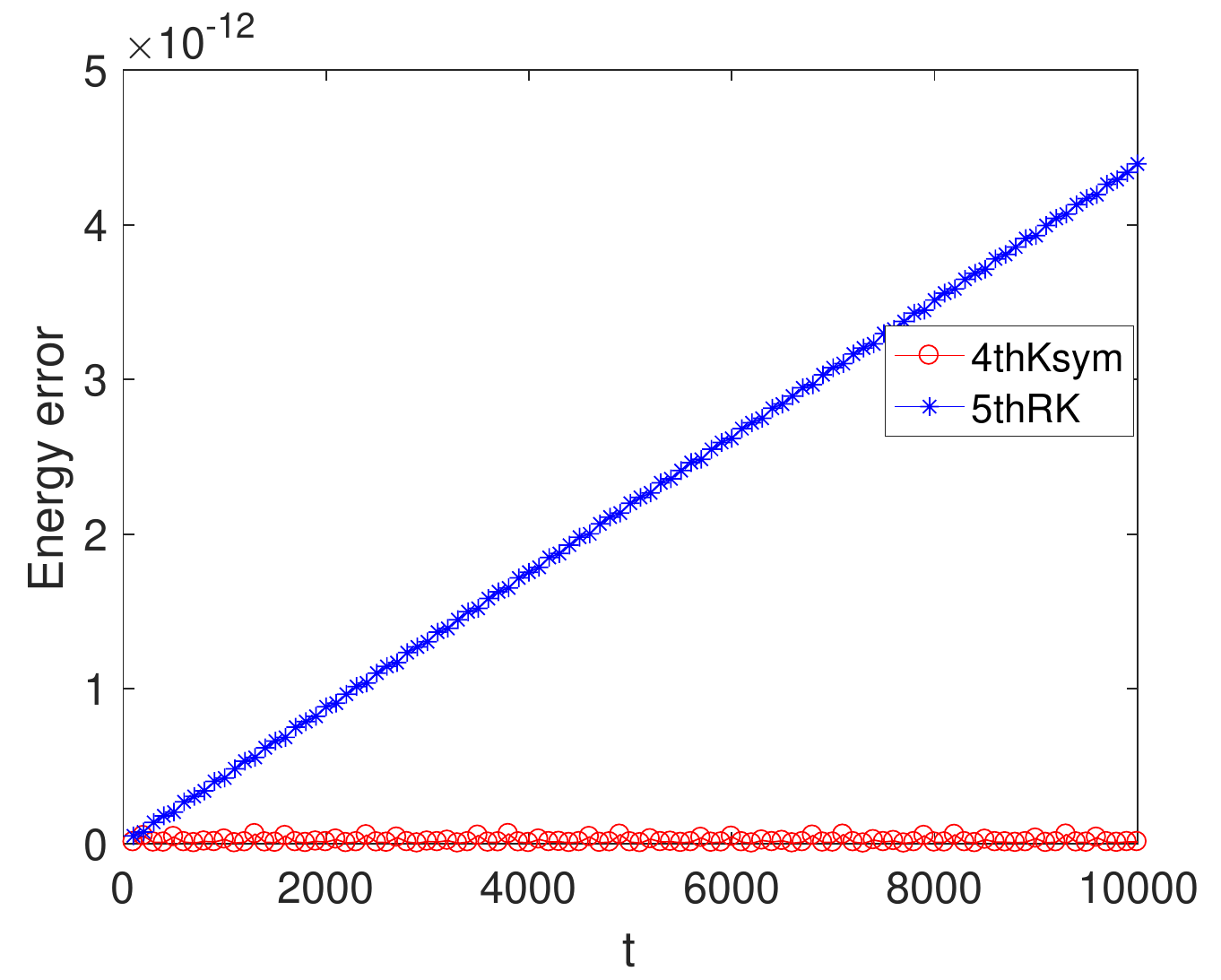}}
\caption{The phase orbit and the energy error obtained by the 2ndKsym method, the 3rdRK method, the 4thKsym method and the 5thRK method. Subfigure (a) and (b) display the orbit projected to $p_1-p_2$ plane obtained by the 2ndKsym method and the 3rdRK method with stepsize $\tau=0.01$, $\Omega=20$ over the time interval $[0,10000]$.  Subfigure (c) and (d) display the relative energy error of the augmented Hamiltonian $\bar{H}$ for the four methods with $\Omega=20$. In subfigure (c), the stepsize $\tau=0.01/4$ while in subfigure (d) $\tau=0.01$. The energy error is represented by $(\bar{H}-\bar{H}_0)/\bar{H}_0$.}
\label{fig:phase2}
\end{figure}

\begin{figure}[p]
\centering
\subfigure[ ]{
\includegraphics[scale=0.42]{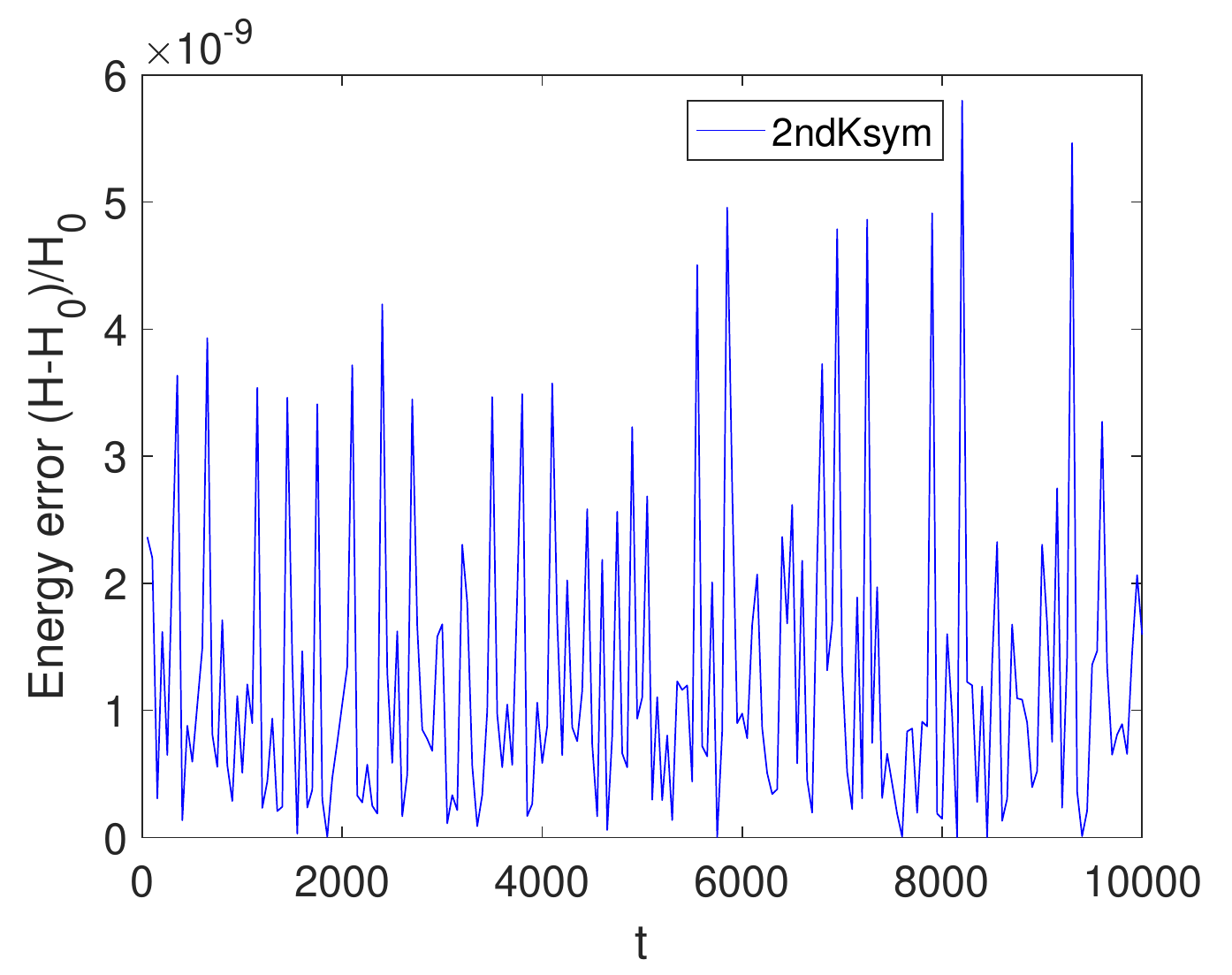}}
\subfigure[ ]{
\includegraphics[scale=0.42]{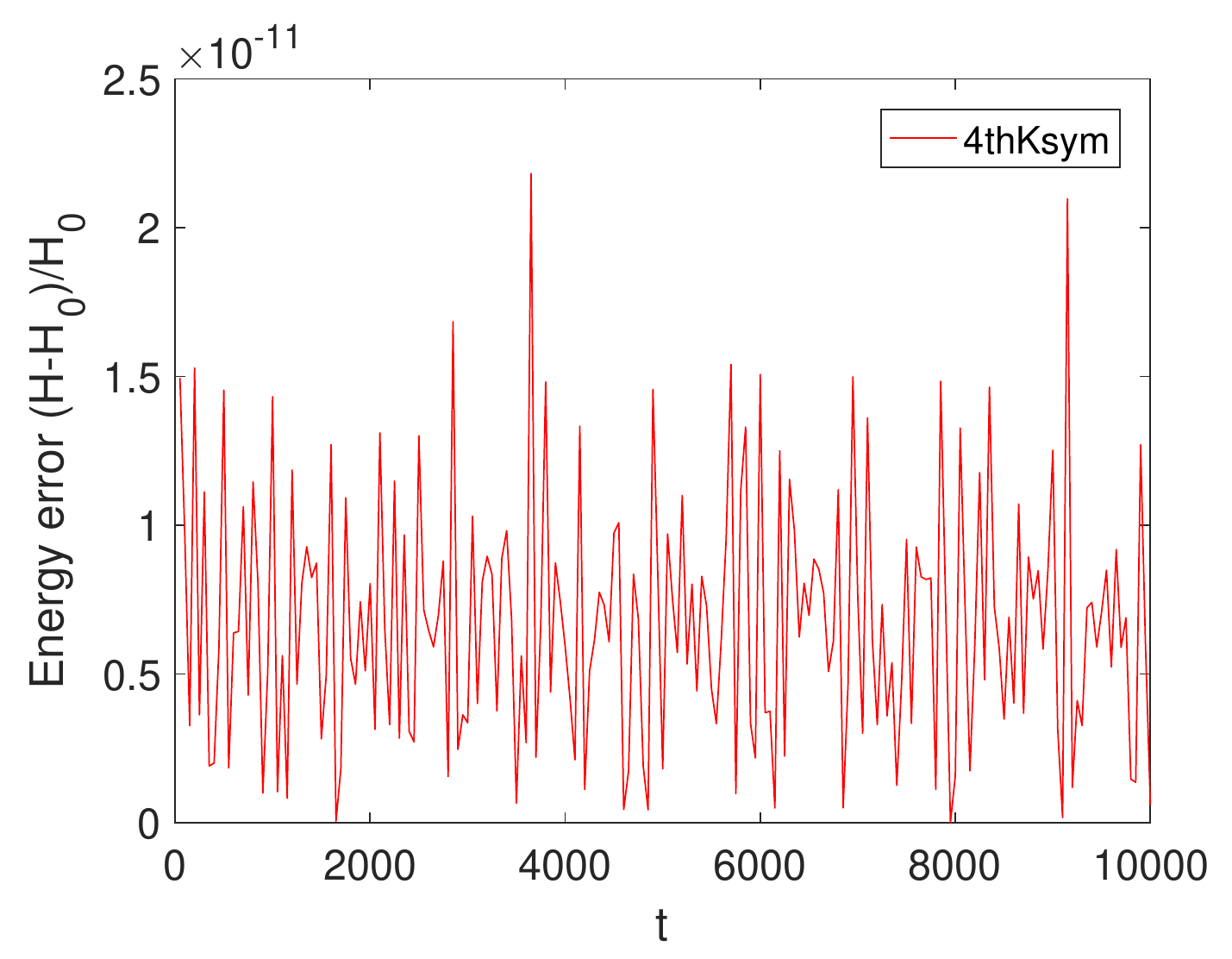}}
\caption{The relative energy error of the original Hamiltonian $H(p_1,p_2,p_3,p_4)$ obtained by the 2ndKsym method and the 4thKsym method. The relative energy error is represented by $(H-H_0)/H_0$. The stepsize is chosen as $\tau=0.01$ and the parameter is $\Omega=20$.}
\label{fig:energy2}
\end{figure}

\begin{figure}[p]
\centering
\subfigure[ ]{
\includegraphics[scale=0.42]{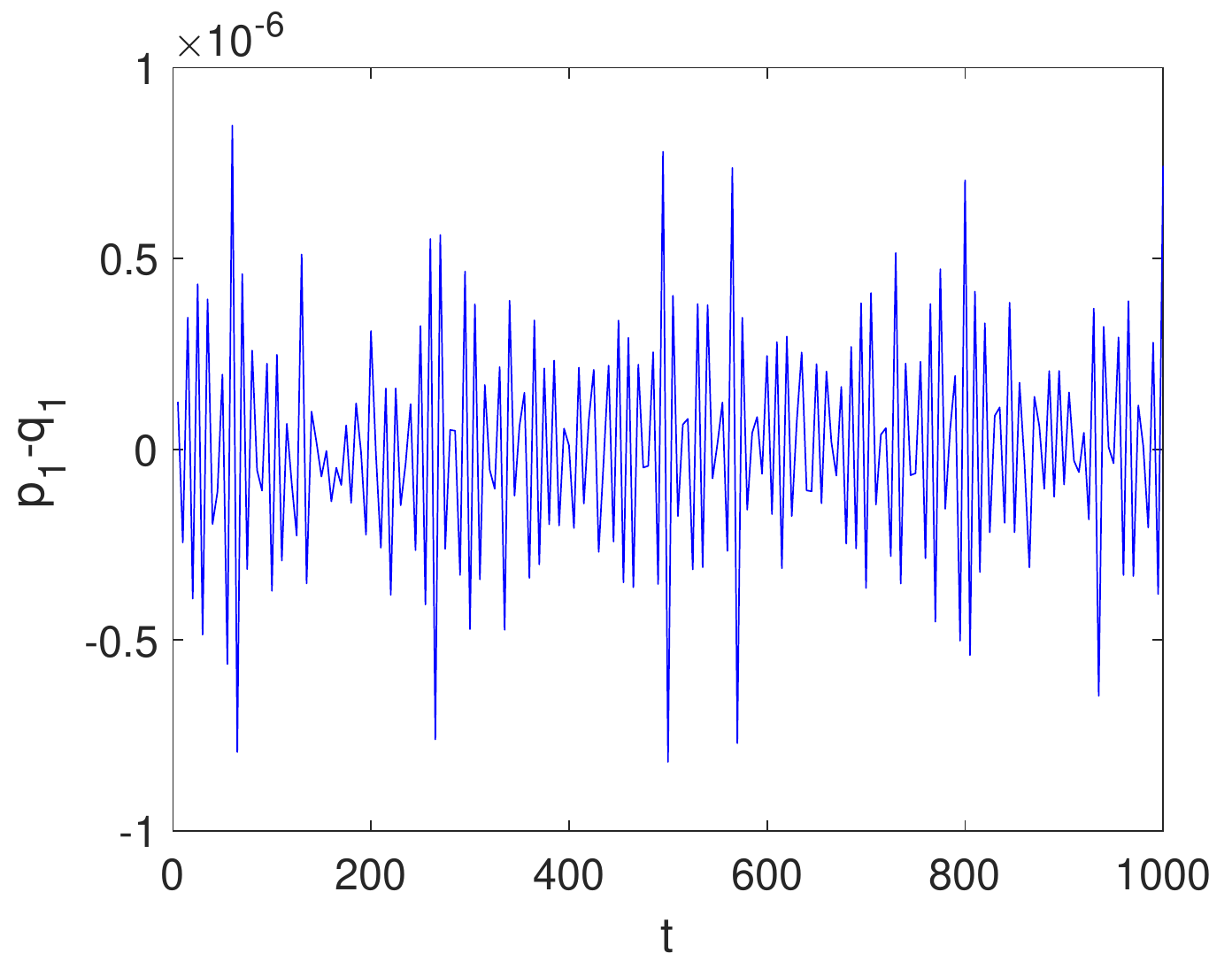}}
\subfigure[ ]{
\includegraphics[scale=0.42]{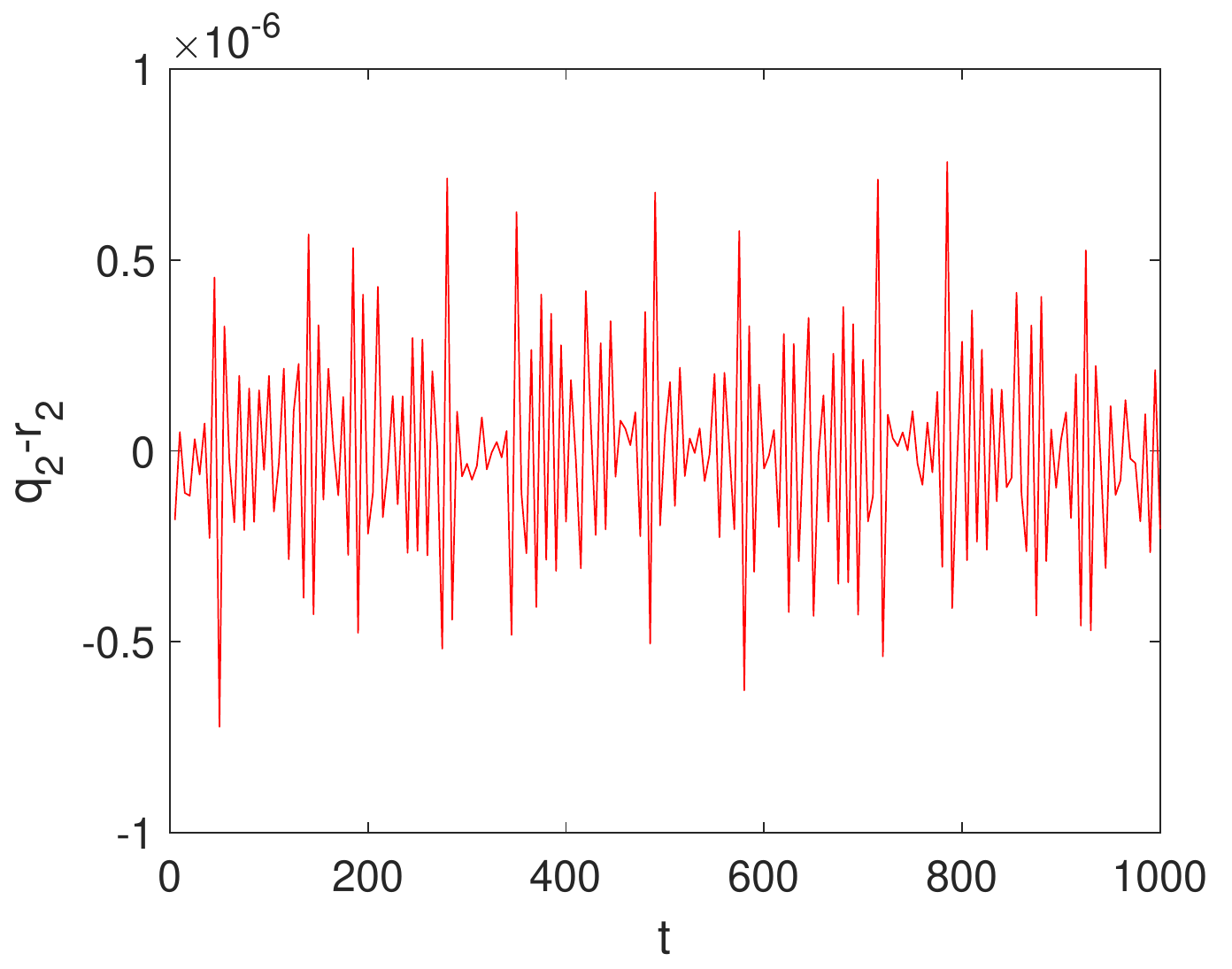}}
\subfigure[ ]{
\includegraphics[scale=0.42]{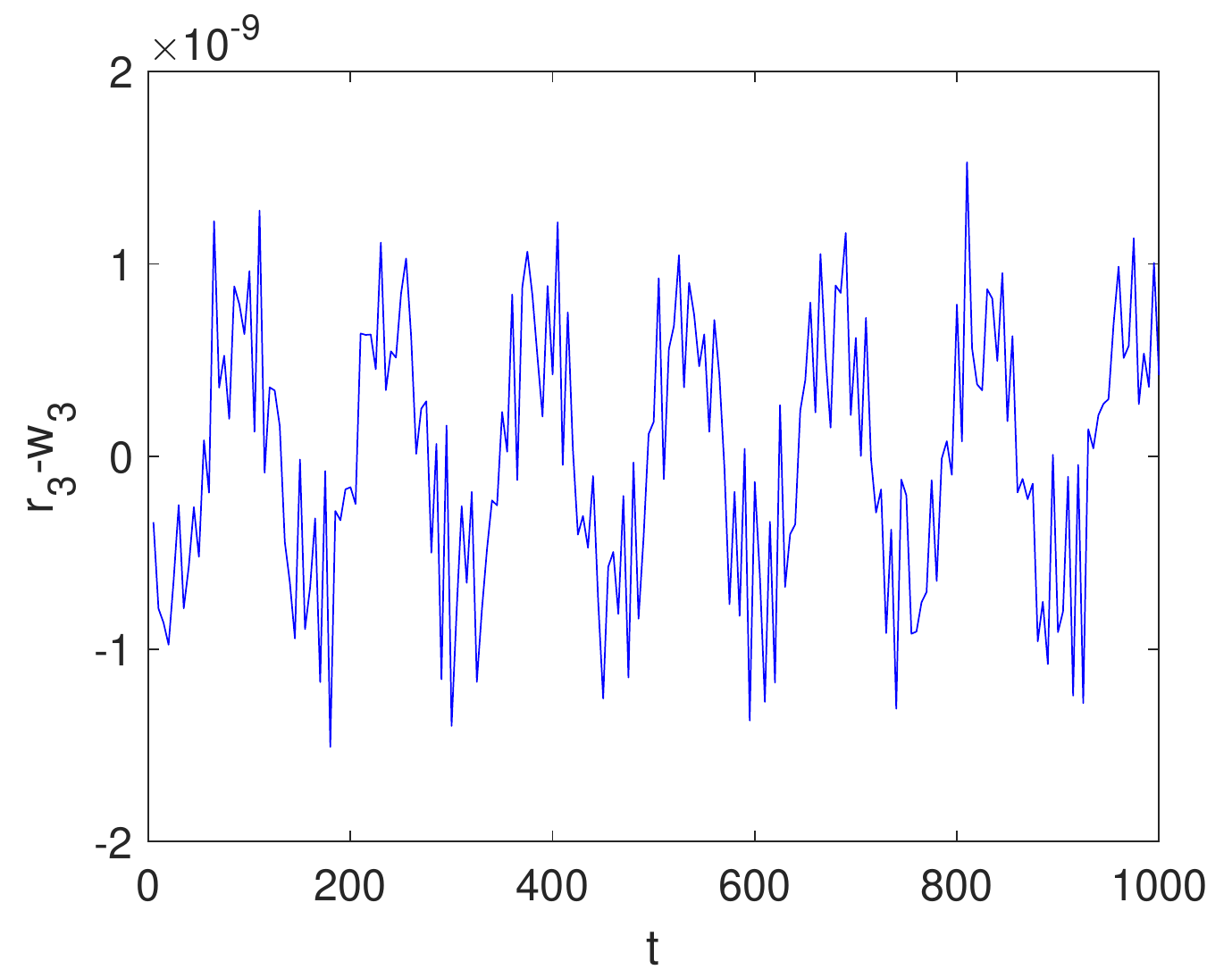}}
\subfigure[ ]{
\includegraphics[scale=0.42]{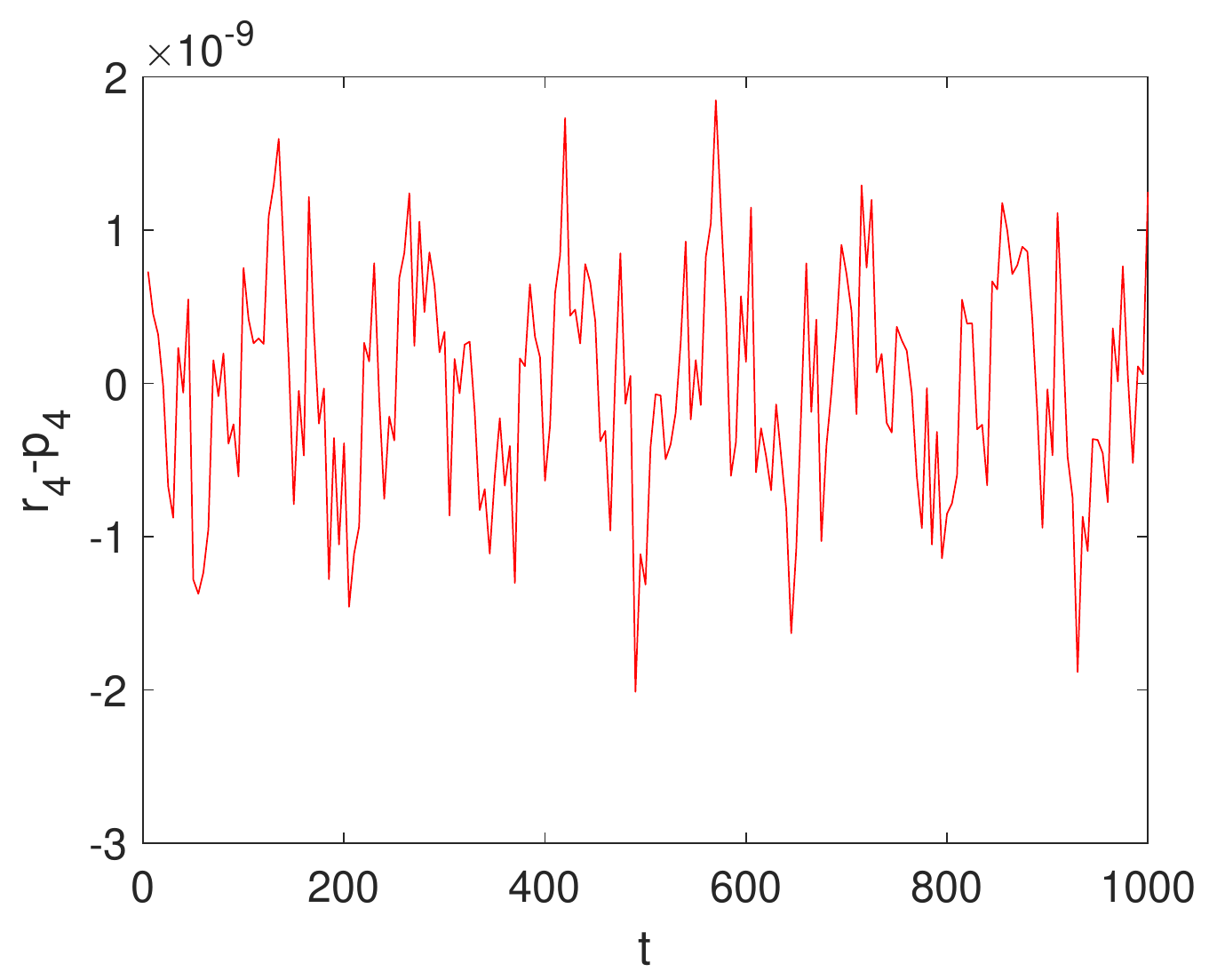}}
\caption{The difference between the four copies for the four variables. Subfigure (a) displays the difference between the first copy and the second copy for the variable $x$ obtained by the 2ndKsym method. Subfigure (b) displays the difference between the second copy and the third copy for $y$ obtained by the 2ndKsym method. Subfigure (c) displays the difference between the third copy and the fourth copy for $z$ obtained by the 4thKsym method, subfigure (b) displays the difference between the first copy and the fourth copy for $u$ obtained by the 4thKsym method. Here the stepsize is $\tau=0.01$ and $\Omega=20$, the final time is $T=1000$.}
\label{fig:pqerror2}
\end{figure}

The initial condition is $x_0=0.003, y_0=0.002, z_0=0.004, u_0=0.005$. The four copies of the original variables $x,y,z,u$ have the same initial condition.  The numerical results obtained by the four numerical methods are displayed in Figure \ref{fig:phase2}-\ref{fig:pqerror2}. The orbits in $p_1-p_2$ plane obtained by the 2ndKsym method and 3rdRK method are displayed in Figure \ref{fig:phase2}. The orbit obtained by the second order K-symplectic method is a closed circle while the orbit obtained by the third order Runge-Kutta method spirals outwards and lose accuracy as can be seen from Figure \ref{fig:phase2}. The relative energy errors of the augmented Hamiltonian $\bar{H}$ obtained by different methods are also shown in Figure \ref{fig:phase1}. The energy errors of K-symplectic methods can be bounded along time while those of the higher order Runge-Kutta methods increase linearly along time. The evolutions of the original $H$ for the two K-symplectic methods are shown in Figure \ref{fig:energy2}. The relative energy errors of the two methods oscillate at vey small amplitudes as can be seen from Figure \ref{fig:energy2}. The differences between the four copies of the original variables $x,y,z,u$ are displayed in Fig. \ref{fig:pqerror2} and they are all bounded at a very small number. Table \ref{table2} displays the CPU time of the four methods under the same stepsize and time interval. It can be seen that the computational cost of the second order explicit K-symplectic method is nearly the same as the third explicit Runge-Kutta method. The CPU time of the fourth order explicit K-symplectic method is three times longer than that of the fifth order explicit Runge-Kutta method. It can be seen from Table \ref{table2} that although we separate the extended system into $d^2$ subsystems, the computational cost is also acceptable. The K-symplectic methods have also shown their superiority in preserving the orbit and the energy over long-term simulation compared with the higher order Runge-Kutta methods.

\begin{table}[htbp]
\begin{small}
\caption{The CPU times of the four methods. The stepsize is $\tau=0.01$ and the time interval is $[0,1000]$.}
\begin{center}
\begin{tabular}{|r|c|c|c|}
\hline
2ndKsym & 3rdRK & 4thKsym & 5thRK\\
\hline
 1.8832 & 1.8378  &  13.2598 & 4.4682 \\
\hline
\end{tabular}
\end{center}
\label{table2}
\end{small}
\end{table}

\section{Conclusion}

We have constructed the explicit K-symplectic methods for non-separable non-canonical Hamiltonian systems. The main technique is extending the phase space to higher dimensional phase space to make the Hamiltonian separable and using the splitting method. If the matrix $K^{-1}$ has some special structure as we mention in Section \ref{Kspecial}, then we make two copies of the phase space. If the matrix $K^{-1}$ does not have special structure, then we make $d$ copies of the phase space where $d$ is the dimension of the phase space. By extending the phase space and imposing some mechanical restraints, then the augmented Hamiltonian becomes separable. We separate the extended system into several subsystems and explicitly solve the subsystems, then the explicit K-symplectic methods can be constructed by composing the exact solution of all subsystems. We have analyzed the situations in which the explicit K-symplectic methods can be constructed.

Explicit K-symplectic methods have been constructed for three non-separable non-canonical Hamiltonian systems and are compared with higher order explicit Runge-Kutta methods. The numerical results show that the K-symplectic methods have superiority in phase orbit tracking and energy conservation over long-term simulation compared with the higher order Runge-Kutta methods. The differences between the several copies of the variables are also bounded by the K-symplectic methods. We have also compared with the CPU times of the K-symplectic methods and the Runge-Kutta methods. The results show that the higher order K-symplectic methods takes more CPU times than the Runge-Kutta methods, but the computational cost is somehow acceptable.

\section*{Acknowledgments}

This research is supported by
the National Natural Science Foundation of China (Grant Nos. 11901564 and 12171466).

\bibliographystyle{abbrv}
\bibliography{main}

\end{document}